\documentstyle[amscd,amssymb,verbatim,diagrams,12pt]{amsart}
\newcommand{\Kr}{\operatorname{Kr}}

\newcommand{\Art}{\operatorname{Art}}
\newcommand{\Sets}{\operatorname{Sets}}

\newcommand{\vareps}{\varepsilon}

\newcommand{\del}{\partial}

\newcommand{\fG}{{\frak G}}

\renewcommand{\mod}{\operatorname{mod}}

\newcommand{\OO}{{\cal O}}
\newcommand{\RR}{{\cal R}}

\newcommand{\Def}{\operatorname{Def}}

\let\liminv\varprojlim

\newcommand{\KK}{{\cal K}}

\newcommand{\be}{{\bf e}}

\newcommand{\G}{{\Bbb G}}

\newcommand{\mg}{{\frak m}}
\newcommand{\hra}{\hookrightarrow}
\newcommand{\lan}{\langle}
\newcommand{\ran}{\rangle}

\newcommand{\CC}{{\cal C}}
\newcommand{\UU}{{\cal U}}
\newcommand{\WW}{{\cal W}}

\newcommand{\Spec}{\operatorname{Spec}}

\newcommand{\Proj}{\operatorname{Proj}}
\renewcommand{\P}{{\Bbb P}}

\newcommand{\si}{\sigma}

\newcommand{\ga}{\gamma}
\newcommand{\de}{\delta}

\newcommand{\im}{\operatorname{im}}

\newcommand{\A}{{\Bbb A}}

\numberwithin{equation}{subsection}

\newcommand{\GL}{\operatorname{GL}}

\newtheorem{thm}{Theorem}[subsection]
\newtheorem{prop}[thm]{Proposition}
\newtheorem{lem}[thm]{Lemma}
\newtheorem{cor}[thm]{Corollary}
{  \theoremstyle{definition}
\newtheorem{defi}[thm]{Definition}
\newtheorem{ex}[thm]{Example}

\newtheorem{rem}[thm]{Remark}
\newtheorem{rems}[thm]{Remarks}
}

\newcommand{\Pf}{\noindent {\it Proof}}
\newcommand{\id}{\operatorname{id}}

\newcommand{\ov}{\overline}

\newcommand{\cusp}{\operatorname{cusp}}
\newcommand{\aff}{\operatorname{aff}}

\newcommand{\Sch}{\operatorname{Sch}}

\renewcommand{\AA}{{\cal A}}

\newcommand{\FF}{{\cal F}}
\newcommand{\EE}{{\cal E}}
\newcommand{\MM}{{\cal M}}
\newcommand{\TT}{{\cal T}}

\newcommand{\HH}{{\cal H}}

\newcommand{\Ext}{\operatorname{Ext}}

\renewcommand{\a}{\alpha}
\renewcommand{\b}{\beta}
\newcommand{\om}{\omega}

\newcommand{\la}{\lambda}

\newcommand{\Z}{{\Bbb Z}}
\newcommand{\Q}{{\Bbb Q}}

\newcommand{\wt}{\widetilde}
\newcommand{\ot}{\otimes}

\newcommand{\sub}{\subset}
\newcommand{\ed}{\qed\vspace{3mm}}

\newcommand{\Qcoh}{\operatorname{Qcoh}}

\newcommand{\cha}{\operatorname{char}}

\newcommand{\ev}{\operatorname{ev}}

\newcommand{\ba}{{\bf a}}

\newcommand{\sslash}{\mathbin{/\mkern-6mu/}}

\title{Moduli of curves with nonspecial divisors and relative moduli of $A_\infty$-structures}
\author{Alexander Polishchuk}
\thanks{Supported in part by NSF grant}

\begin{document}
\begin{abstract} In this paper for each $n\ge g\ge 0$ we consider the moduli stack $\wt{\UU}^{ns}_{g,n}$
of curves $(C,p_1,\ldots,p_n,v_1,\ldots,v_n)$ of
arithmetic genus $g$ with $n$ smooth marked points $p_i$ and nonzero tangent vectors $v_i$ at them, such that
the divisor $p_1+\ldots+p_n$ is nonspecial (has no $h^1$) and ample. With some mild restrictions on the characteristic 
we show that it is a scheme, affine over the Grassmannian $G(n-g,n)$. We also construct an isomorphism of
$\wt{\UU}^{ns}_{g,n}$ with a certain relative moduli of $A_\infty$-structures (up to an equivalence) over a family of
graded associative algebras parametrized by $G(n-g,n)$.
\end{abstract}

\maketitle

\section*{Introduction}

This paper continues the study of connections between moduli spaces of curves and $A_\infty$-algebras,
started in \cite{P-ainf}. Recall that in \cite{P-ainf} we gave an interpretation of the moduli of curves
of arithmetic genus $g$ with $g$ (distinct smooth) marked points forming a nonspecial divisor, as a certain moduli
space of $A_\infty$-algebras. In the present paper we consider a generalization of this picture to the case of
curves with $n$ marked points, where $n\ge g$. 

The second motivation for this work is the relation, pointed out in \cite{LP}, between the 
moduli space of curves $(C,p_1,\ldots,p_n)$ of arithmetic genus $1$ such that $H^1(C,\OO(p_i))=0$ for every $i$
and one of the moduli spaces studied by Smyth in \cite{Smyth-I} and \cite{Smyth-II}. For each $1\le m<n$
he constructed an alternate compactification $\ov{\MM}_{1,n}(m)$
of the moduli space of $n$-pointed curves of genus $1$ consisting of {\it $m$-stable} curves. The moduli space
that shows up in \cite{LP} is $\ov{\MM}_{1,n}(n-1)$. In the present paper we construct a bigger moduli stack of curves
which should contain open substacks closely related to Smyth's moduli spaces for all $m\ge (n-1)/2$
(see Section \ref{g=1-any-n-sec}).

Let us fix $n\ge g$. The main object of study of this paper is the moduli stack
$\UU^{ns}_{g,n}$ of $(C,p_1,\ldots,p_n)$ such that $H^1(C,\OO(p_1+\ldots+p_n))=0$ and $\OO(p_1+\ldots+p_n)$ is ample,
and the $\G_m^n$-torsor over it, $\wt{\UU}^{ns}_{g,n}$, corresponding to choices of nonzero
tangent vectors $v_1,\ldots,v_n$ at the marked points. 

Note that the vanishing of $H^1(C,\OO(p_1+\ldots+p_n))$ is equivalent to the surjectivity of the map 
\begin{equation}\label{connecting-hom-eq}
H^0(C,\OO_C(p_1+\ldots+p_n)/\OO_C)\to H^1(C,\OO_C).
\end{equation}
Hence, its kernel is $n-g$-dimensional. Thus, we have a natural morphism
\begin{equation}\label{Gr-map}
\pi:\wt{\UU}_{g,n}^{ns}\to G(n-g,n),
\end{equation}
where $G(n-g,n)$ is the Grassmannian of $(n-g)$-dimensional subspaces in the $n$-dimensional space,
associating with $(C,p_1,\ldots,p_n,v_1,\ldots,v_n)$ the kernel of the map \eqref{connecting-hom-eq},
where $H^0(C,\OO_C(p_1+\ldots+p_n)/\OO_C)$ is trivialized using the basis $v_1,\ldots,v_n$.

Note that some closely related moduli stacks were considered in \cite{P-krich}. Namely, 
the preimages of the standard cells in $G(n-g,n)$ under $\pi$ are the open 
substacks $\wt{\UU}^{ns}_{g,n}(S)\sub \wt{\UU}^{ns}_{g,n}$, 
for subsets $S\sub\{1,\ldots,n\}$ such that $|S|=g$,
given by the condition $H^1(C,\OO_C(\sum_{i\in S}p_i))=0$.
These stacks are precisely the stacks $\wt{\UU}^{ns}_{g,n}(\ba)$, defined in \cite{P-krich}
for collections $\ba=(a_1,\ldots,a_n)\in \Z_{\ge 0}^n$ such that $\sum_i a_i=g$, in the case when
each $a_i$ is either $0$ or $1$.

Working over $\Q$ we proved in \cite{P-krich} that each $\wt{\UU}^{ns}_{g,n}(S)$ is in fact an affine scheme of
finite type, and identified it with the quotient of a certain locally closed subset of the Sato Grassmannian
of subspaces in $\HH=\bigoplus_{i=1}^n k((t_i))$ by the free action of the group of changes of variables.
The first result of this paper, Theorem A below, gives analogous statements for $\wt{\UU}^{ns}_{g,n}$
(in this case the morphism $\pi$ is affine of finite type).

As in \cite{P-krich}
we consider the closed subset $ASG$ of the Sato Grassmannian consisting of $W$ that are subalgebras of $W$.
Let $ASG^{ns}\sub ASG$ be the open subset consisting of $W$ such that
$$W\cap \bigoplus_{i=1}^n k[[t_i]]=k, \ \dim(\HH/(W+\bigoplus_{i=1}^n k[[t_i]]))=g, \text{ and }
\ \HH=W+\bigoplus_{i=1}^n t_i^{-1}k[[t_i]].
$$
There is a natural action on $ASG^{ns}$ of the group $\fG$ of changes of variables of the form 
$t_i\mapsto t_i+c_{1i}t_i^2+c_{2i}t_i^3+\ldots$, $i=1,\ldots,n$.

\medskip

\noindent
{\bf Theorem A (see Theorem \ref{big-moduli-thm}).} {\it (i) Assume that either 
\begin{itemize}
\item $n\ge g\ge 1$, $n\ge 2$ and the base is $\Spec(\Z[1/2])$, or 
\item $n=g=1$ and the base is $\Spec(\Z[1/6])$, or 
\item $g=0$, $n\ge 2$ and the base is $\Spec(\Z)$.
\end{itemize}
Then the stack $\wt{\UU}^{ns}_{g,n}$ is isomorphic to a scheme, affine of finite type over the Grassmannian $G(n-g,n)$, so that the preimages of the standard open cells $U_S\sub G(n-g,n)$,
for $S\sub\{1,\ldots,n\}$, $|S|=g$, are the moduli schemes
$\wt{\UU}_{g,n}(\ba_S)$, where $\ba_S$ has $1$'s at the places corresponding to $S$. 

\noindent
(ii) Now let us work over $\Spec(\Q)$. Then the action of $\fG$ on 
$ASG^{ns}$ is free, and the Krichever map induces an isomorphism
$$\wt{\UU}^{ns}_{g,n}\simeq ASG^{ns}/\fG.$$
}

Note that part (i) of this Theorem
is an improvement of \cite[Thm.\ 1.2.4]{P-ainf} (where the case $n=g$ was considered, but over $\Z[1/6]$)
and of the special case of \cite[Thm.\ A(i)]{P-krich} when $\ba$ is a collection of $0$'s and $1$'s.

Next, generalizing the work \cite{P-ainf} (corresponding to the case $n=g$), we consider $A_\infty$-algebras
associated with curves $(C,p_\bullet,v_\bullet)$ in $\wt{\UU}^{ns}_{g,n}$.
Namely, we consider the object
\begin{equation}\label{G-generator-eq}
G=\OO_C\oplus\OO_{p_1}\oplus\ldots\oplus\OO_{p_n}
\end{equation}
in the perfect derived category of $C$ and consider the natural minimal $A_\infty$-structure on the corresponding
algebra $\Ext^*(G,G)$ (which arises from a dg-model of this $\Ext$-algebra and is 
defined uniquely up to a gauge equivalence).
The key observation is that the associative algebra structure on $\Ext^*(G,G)$ depends only on the corresponding
$(n-g)$-dimensional subspace in $k^n$. 

More precisely, let $Q_n$ be the quiver with $n+1$ vertices marked as $\OO,\OO_{p_1},\ldots,\OO_{p_n}$ and with the arrows 
$$A_i:\OO_{p_i}\to \OO, \ \ B_i:\OO\to \OO_{p_i}, \ i=1,\ldots,n.$$
Let $J_0$ be the two-sided ideal in the path algebra $k[Q_n]$ of $Q_n$ generated by the elements
$$A_iB_iA_i, B_iA_iB_i, A_iB_j,$$
where $i\neq j$.
For an $(n-g)$-dimensional subspace $W\sub k^n$
we define $J_W\sub k[Q_n]$ to be  the ideal generated by $J_0$ together with the additional relations
$$\sum x_i B_iA_i=0 \ \ \text{ for every } \sum x_i e_i\in W,$$
and consider the corresponding quotient algebra
\begin{equation}\label{E-W-eq}
E_W=k[Q_n]/J_W
\end{equation}
We equip $E_W$ with the $\Z$-grading by $\deg(A_i)=0$, $\deg(B_i)=1$.

Now for a curve $(C,p_\bullet,v_\bullet)\in \wt{\UU}^{ns}_{g,n}$ there is a canonical
isomorphism of associative algebras 
$$\Ext^*(G,G)\simeq E_W,$$
for $W=\pi(C,p_\bullet,v_\bullet)$. Thus, from such a curve we get an $A_\infty$-structure on
the algebra $E_W$. 

The family of associative algebras $E_W$ defines a sheaf of $\OO$-algebras $\EE_{g,n}$
over $G(n-g,n)$. Extending the techniques developed in \cite{P-ainf} we consider the relative moduli space 
$\MM_\infty$ over $G(n-g,n)$, 
classifying minimal $A_\infty$-structures on the fibers of $\EE_{g,n}$ (for a precise definition see 
Def.\ \ref{ainf-moduli-def}).
We prove that in fact $\MM_\infty$ is an affine scheme over $G(n-g,n)$ (over $\Z[1/6]$), and the above
construction of $A_\infty$-structures associated with curves gives an isomorphism of the moduli spaces.

\noindent
{\bf Theorem B} (see Theorem \ref{ainf-thm}). {\it Under the assumptions of Theorem A(i) we have a natural isomorphism
$$\wt{\UU}^{ns}_{g,n}\rTo{\sim} \MM_\infty$$
of affine schemes over $G(n-g,n)$, compatible with the $\G_m^n$-action, where $(\la_i)\in\G_m^n$ acts on 
$\wt{\UU}^{ns}_{g,n}$ by rescaling the tangent vectors at the marked points and on $\MM_\infty$ by the rescalings
$$A_i\mapsto A_i, \ B_i\mapsto \la_iB_i.$$
}

Note that in the case $n=g\ge 2$ (resp., $n\ge 3$, $g=0$)
this result is a strengthening of a similar isomorphism in \cite[Thm.\ A]{P-ainf} (resp.,
\cite[Thm.\ 5.2.1]{P-ainf}), since we now work over $\Z[1/2]$ (resp., $\Z$), not over a field.

The construction of the scheme $\MM_\infty$ of minimal $A_\infty$-structures on $\EE_{g,n}$ is a particular case of a
more general construction (see Theorem \ref{gen-ainf-thm})
of the affine scheme classifying equivalence classes of minimal $A_\infty$-structures on a sheaf of
$\OO$-algebras $\EE$ over a scheme $S$, such that $\EE$ is locally free of finite rank as an $\OO$-module,
and the associative algebras $\EE_s$ for $s\in S$ satisfy the following vanishing condition:
\begin{equation}\label{HH-vanishing-eq}
HH^i(\EE_s)_{<0}=0 \ \text{ for } i\le 1.
\end{equation}

As in \cite{P-ainf} the important part of the proof of Theorem B is identifying the curves in $\wt{\UU}^{ns}_{g,n}$
such that the corresponding $A_\infty$-algebras are homotopically trivial. In the case $n=g$, considered in \cite{P-ainf},
there is only one such curve, $C_g^{\cusp}$, which is the union of $g$ cuspidal curves of genus $1$, glued at the cusp. 
In general, there is a family of such curves, parametrized by $G(n-g,n)$. Namely, these 
are precisely the invariant points of the action of the diagonal $\G_m\sub\G_m^n$
on $\wt{\UU}^{ns}_{g,n}$.
We refer to curves in this family as {\it special curves}. 
In Theorem \ref{big-moduli-thm} we prove that special curves form a section of the projection \eqref{Gr-map}.
Special curves are used in the proof of Theorem B as follows. Since the $\G_m$-action contracts both
spaces, $\wt{\UU}^{ns}_{g,n}$ and $\MM_\infty$, to the $\G_m$-invariant locus (which is $G(n-g,n)$ in both cases),
it is enough to study deformations of each special curve and show that they precisely correspond to deformations of
the $E_W$ as an $A_\infty$-algebra. This is done using the same ideas as in the case $n=g$, although there are
some new features that appear because we now work with a family of associative algebras (see Proposition
\ref{deform-prop}). 

Note that our moduli scheme of $A_\infty$-structures $\MM_\infty$ has a natural extension to a derived stack
that can be constructed as in \cite[(3.2)]{CK}. It would be interesting to find an interpretation of this derived extension 
in terms of moduli of curves.

The paper is organized as follows. Section \ref{curves-sec} is devoted to geometric aspects of the moduli
stacks $\wt{\UU}^{ns}_{g,n}$. In particular, in Section \ref{special-curves-sec} we describe a family of special
curves in $\wt{\UU}^{ns}_{g,n}$. Then in Section \ref{constr-sec} we prove Theorem A.
In Section \ref{gluing-sec} we describe the natural gluing morphism that associates with a pair of curves
from the moduli spaces $\wt{\UU}^{ns}_{g_1,n_1}$ and $\wt{\UU}^{ns}_{g_2,n_2}$, each equipped with an additional
point different from all the marked points, a glued curve in $\wt{\UU}^{ns}_{g_1+g_2,n_1+n_2}$.
In Sections \ref{g=1-n=2-sec} and \ref{g=1-any-n-sec} we study the case $g=1$: we describe explicitly
the space $\wt{\UU}^{ns}_{1,2}$, as well as construct regular morphisms from the Smyth's moduli spaces
of $m$-stable curves to $\wt{\UU}^{ns}_{1,n}$ for $m\ge (n-1)/2$.
Section \ref{ainf-sec} is devoted to the relative moduli of $A_\infty$-structures. After proving some technical results
in Section \ref{nice-quotients-sec}, we give  in Section \ref{gen-ainf-moduli-sec}
a general construction of the affine
scheme parametrizing $A_\infty$-structures over a given family of associative algebras
(under the assumption \eqref{HH-vanishing-eq}). Finally, in Section \ref{ainf-curve-sec} we prove Theorem B.


\section{Moduli of curves with nonspecial divisors}\label{curves-sec}

\subsection{Some special curves}\label{special-curves-sec}

First, we are going to construct some special curves that will play an important role in our study of the moduli spaces
$\wt{\UU}^{ns}_{g,n}$.

\begin{defi}\label{sp-curve-def} (i) Let $x_1,\ldots,x_n$ be independent variables, and let $R$ be a commutative ring. 
For a subset $S\sub\{1,\ldots,n\}$ with $|S|=g$
let us consider the subalgebra in
$\bigoplus_{i=1}^nR[x_i]$ given by
$$B(S):=R\cdot 1+\bigoplus_{i\in S} x_i^2k[x_i].$$
Next, let $(\ov{h}_j)_{1\le j\le n, j\not\in S}$ be a collection of linear forms in $(x_i)_{i\in S}$ with coefficients in $R$.
We define the $R$-algebra $A(\ov{h}_\bullet)$ as the $B(S)$-subalgebra in $\bigoplus_{i=1}^n R[x_i]$ 
generated by the elements $h_j=x_j+\ov{h}_j$, $j\not\in S$.
We view $A(\ov{h}_\bullet)$ as a graded $R$-algebra, where $\deg(x_i)=1$.

\noindent
(ii) We define two curves, one affine and another projective over $R$, by 
$$C^{\aff}(\ov{h}_\bullet)=\Spec(A(\ov{h}_\bullet)),$$ 
$$C(\ov{h}_\bullet)=\Proj(\RR(A(\ov{h}_\bullet))),$$ 
where $\RR(A(\ov{h}_\bullet))=\bigoplus_{m\ge 0} F_m$ is the Rees algebra associated
with the increasing filtration $(F_m)$ on $A(\ov{h}_\bullet)$ coming from the grading.
Note that $C^{\aff}(\ov{h}_\bullet)$ is an affine open in $C(\ov{h}_\bullet)$.
we have an action of $\G_m$ on these curves associated with the grading on $A(\ov{h}_\bullet)$.
\end{defi}


\begin{prop}\label{special-curves-prop} 
(i) For any matrix $(a_{ij})_{i\in S,j\not\in S}$ with entries in $R$, let $A$ be the graded algebra 
defined by the generators
$(f_i,h_i,h_{S,j})_{i\in S,j\not\in S}$ subject to the equations 
\begin{align}
\begin{split}
\label{Gm-inv-curves-eq}
&f_if_{i'}=0, \ \ f_ih_{i'}=0, \ \ h_ih_{i'}=0, \ \ h_i^2=f_i^3,\\
&h_{S,j}h_{S,j'}=\sum_{i\in S}a_{ij}a_{ij'}f_i,\\
&f_ih_{S,j}=a_{ij}h_i,\\
&h_ih_{S,j}=a_{ij}f_i^2,
\end{split}
\end{align}
where $i,i'\in S$, $i\neq i'$, $j,j'\not\in S$, $j\neq j'$, and the grading is given by
$$\deg(h_{S,j})=1, \ \deg(f_i)=2, \ \deg(h_i)=3.$$
Then there is an injective homomorphism of graded $R$-algebras 
$$\rho:A\to \bigoplus_{i=1}^n R[x_i],$$
such that
$$\rho(f_i)=x_i^2, \ \ \rho(h_i)=x_i^3, \ i\in S,$$
$$\rho(h_{S,j})=x_j+\sum_{i\in S} a_{ij}x_i, \ j\not\in S,$$
inducing an isomorphism $A\simeq A(\ov{h}_\bullet)$, where
$$\ov{h}_j=\sum_{i\in S} a_{ij} x_i.$$
The elements 
\begin{equation}\label{Grobner-basis-A-h-eq}
(f_i^n, f_i^nh_i, h_{S,j}^m), \ \ i\in S, j\not\in S, m\ge 1, n\ge 0,
\end{equation}
form a basis of $A$ over $R$.

\noindent
(ii) Assume now $R=k$, where $k$ is a field. Let $C=C(\ov{h}_\bullet)$, $C^{\aff}=C^{\aff}(\ov{h}_\bullet)$ be as in
Definition \ref{sp-curve-def}.
Then $C\setminus C^{\aff}$ consists of $n$ smooth points $p_1,\ldots,p_n$ on $C$.
Furthermore, $C$ is a union of $n$ components $C_i$, 
joined in a single point $q$, which is the only singular point
of $C$ (with $p_i\in C_i\setminus \{q\}$). Each component $C_i$ is either $\P^1$, or the cuspidal curve of arithmetic
genus $1$.
\end{prop}

\Pf . (i) It is easy to see that $\rho$ is well defined, and that the elements \eqref{Grobner-basis-A-h-eq}
span $A$ over $R$. On the other hand, one immediately checks that their images under $\rho$
are linearly independent over $R$ and generate the subalgebra $A(\ov{h}_\bullet)$. This implies our
assertions.

\noindent
(ii) The complement $C\setminus C^{\aff}$ is naturally identified with 
$$\Proj A(\ov{h}_\bullet)\simeq \Proj(\bigoplus_{i=1}^n k[x_i])=\sqcup_{i=1}^n p_i.$$
Let us consider the projection from $C^{\aff}(\ov{h}_\bullet)$ to the union of $g$ cuspidal 
curves $h_i^2=f_i^3$ (given by the projection to coordinates $(f_i,h_i)$).
The fiber over the cusp is the union of the coordinate axes in $\A^{n-g}$. We claim that over the complement to the cusp
this projection is an isomorphism.
Say, $h_i=f_i=0$ for all $i\neq i_0$, and $h_{i_0}$, $f_{i_0}$ are invertible. 
Then the equations
$$h_{S,j}h_{S,l}=a_{i_0,j}a_{i_0,l}f_{i_0}, \ f_{i_0}h_{S,j}=a_{i_0,j}h_{i_0}, \ h_{i_0}h_{S,j}=a_{i_0,j}f_{i_0}^2$$
are equivalent to
$$h_{S,j}=a_{i_0,j}\frac{h_{i_0}}{f_{i_0}},$$
for $j\not\in S$.

Thus, $C^{\aff}$ has $n$ irreducible rational components $C^{\aff}_i$, all joined at one point, where $h_i=f_i=h_{S,j}=0$.
More precisely, for $j\not\in S$, we have $C^{\aff}_j\simeq\A^1$: $C^{\aff}_j$ is given by the equations
$f_i=h_i=0$, $h_{S,l}=0$ for $l\neq j$, and $x_j=h_{S,j}$ is the coordinate on it. 

For $i\in S$ there are two cases:

\noindent
Case 1. There exists $j\not\in S$ such that $a_{i,j}\neq 0$. Then $C^{\aff}_i\simeq \A^1$ with the coordinate 
$x_i=h_{S,j}/a_{i,j}$. Note that if $a_{i,l}\neq 0$ for some other $l\not\in S$
then $h_{S,j}/a_{i,j}=h_{S,l}/a_{il}$ on $C^{\aff}_i$.

\noindent
Case 2. $a_{ij}=0$ for all $j\not\in S$. Then all $h_{S,j}=0$ on $C^{\aff}_i$, so $C^{\aff}_i$ is the cuspidal curve
with the coordinate $x_i$ on the normalization of $C^{\aff}_i$ such that $x_i^2=f_i$, $x_i^3=h_i$.

Note the map $\rho$ is precisely the map associating with a function $f\in A$ its 
pull-backs to normalizations of the irreducible components $C^{\aff}_1,\ldots,C^{\aff}_n$, where
we choose coordinates $x_i$ as above. Let $C_i$ be the closure of $C^{\aff}_i$.
The point at infinity $p_i\in C_i$ corresponds 
either to the point on $\P^1$ where $x_i=\infty$, or to the infinite point on the projective closure of the cuspidal curve
(corresponding to $x_i=\infty$ on its normalization). In particular, all $p_i$ are smooth.
\ed

Note that the curve $C(\ov{h})$ is determined by the corresponding subspace 
$$W:=\lan x_j+\ov{h}_j \ |\ j\not\in S\ran\sub \bigoplus_{i=1}^n k\cdot x_i\simeq k^n,$$
which can be any point of the open cell in $G(n-g,n)$ where $W+k^S=k^n$. Later we will show that each
curve $C(\ov{h})$ with the marked points $p_i$ and the natural tangent vectors at them induced by $x_i^{-1}$ (viewed
as elements of the local ring of $C(\ov{h})$ at $p_i$) defines a point in $\wt{\UU}_{g,n}^{ns}$, which we will denote
simply by $C_W$. We will refer to the curves of the form $C_W$ as {\it special curves}.

\begin{rem} The special curves above are particular cases of the curves considered in \cite[Sec.\ 2.1]{P-krich}.
Note that in the case $g=0$ there is a unique special curve for each $n$: the union of $n$ projective lines joined
in a rational $n$-fold point.
\end{rem}

\subsection{Moduli spaces}\label{constr-sec}

Recall that for each collection $a_1,\ldots,a_n\ge 0$ such that $a_1+\ldots+a_n=g$, we
considered in \cite{P-krich} the stack $\wt{\UU}_{g,n}^{ns}(a_1,\ldots,a_n)$ of curves $(C,p_1,\ldots,p_n)$
of arithmetic genus $g$ such that $H^1(C,\OO(\sum a_ip_i))=0$ and $\OO(p_1+\ldots+p_n)$ is ample, equipped
with nonzero tangent vectors at the marked point.
Working over $\Z[1/N]$ for sufficiently divisible $N$, 
we proved that all of these are affine schemes of finite type, and related them (working over $\Q$) to 
certain subschemes of the Sato Grassmannian via the Krichever map.

The moduli spaces considered below are glued from various $\wt{\UU}_{g,n}^{ns}(a_1,\ldots,a_n)$, where $a_i$'s
are all $0$'s and $1$'s.


\begin{defi}
Let us denote by $\UU^{ns}_{g,n}$ the moduli stack of (reduced connected projective)
curves $(C,p_1,\ldots,p_n)$ of arithmetic genus $g$ and $n$ smooth distinct marked points, such that
$H^1(C,\OO(p_1+\ldots+p_n))=0$ and $\OO_C(p_1+\ldots+p_n)$ is ample.
Let $\wt{\UU}^{ns}_{g,n}$ denote the $\G_m^n$-torsor over $\UU^{ns}_{g,n}$
corresponding to choices of nonzero tangent vectors 
at all the marked points.
\end{defi}

For each subset $S\sub\{1,\ldots,n\}$ such that $|S|=g$, we
consider the open substack $\wt{\UU}_{g,n}(S)\sub\wt{\UU}_{g,n}^{ns}$ corresponding
to curves $(C,p_1,\ldots,p_n)$ for which $H^1(C,\OO_C(\sum_{i\in S}p_i))=0$. 
This is equivalent to requiring that $H^0(C,\OO_C(\sum_{i\in S}p_i)/\OO_C)$ surjects onto $H^1(C,\OO_C)$.
Thus, we have 
$$\wt{\UU}_{g,n}(S)=\pi^{-1}(U_S)$$
where $U_S\sub G(n-g,n)$ is the open cell in
the Grassmannian corresponding to subspaces $W\sub k^n$ such that $W+k^S=k^n$.

It follows that $\wt{\UU}_{g,n}^{ns}$ is the union of open substacks
$\wt{\UU}_{g,n}(S)$, where $S$ ranges over subsets of $\{1,\ldots,n\}$
such that $|S|=g$.
On the other other hand, by definition, we have 
$$\wt{\UU}_{g,n}(S)=\wt{\UU}^{ns}_{g,n}(a_1,\ldots,a_n),$$
where $a_i=1$ for $i\in S$ and $a_j=0$ for $j\not\in S$.

We have a natural action of $\G_m^n$ on $\wt{\UU}_{g,n}^{ns}$, so that for
$\la=(\la_1,\ldots,\la_n)$,
$$\la (C,p_1,\ldots,p_n,v_1,\ldots,v_n)=(C,p_1,\ldots,p_n\la_1^{-1}v_1,\ldots,\la_n^{-1}v_n).$$
Note that the \eqref{Gr-map} is $\G_m^n$-equivariant, where $\G_m^n$ acts on $G(n-g,n)$ via the embedding 
$\G_m^n\sub\GL(n)$
and the natural action of $\GL(n)$ on the Grassmannian.

As in \cite{P-krich}, we consider the locally closed subset $SG_1(g)$ in the Sato Grassmannian of 
subspaces of $\HH=\bigoplus_{i=1}^n k((t_i))$, consisting of $W$ such that $1\in W$,
\begin{itemize}
\item $W\cap \HH_{\ge 0}=\lan 1\ran$ and
\item $\dim \HH/(W+\HH_{\ge 0})=g$,
\end{itemize}
where $\HH_{\ge 0}=\bigoplus_{i=1}^n k[[t_i]]$.
We denote by $ASG\sub SG_1(g)$ the closed subscheme
consisting of $W$ which are subalgebras in $\HH$.
We denote by $SG^{ns}\sub SG_1(g)$ the open subset consisting of $W$ such that 
$$W+\bigoplus_{i=1}^n t_i^{-1}k[[t_i]]=\HH,$$
and we set $ASG^{ns}=ASG\cap SG^{ns}$.
All of these schemes can be defined over $\Z$ (see \cite[Sec.\ 1.1]{P-krich} for details).

Let $\UU^{ns,(\infty)}_{g,n}$ be the torsor over $\UU^{ns}_{g,n}$ corresponding to choices of formal
parameters $(t_1,\ldots,t_n)$ at the marked points.
Using \cite[Prop.\ 1.1.5]{P-krich} we see that there is a natural morphism ({\it Krichever map})
\begin{equation}\label{Krich-map}
\Kr:\UU^{ns,(\infty)}_{g,n}\to ASG^{ns}: (C, p_\bullet, t_\bullet)\mapsto H^0(C\setminus\{p_1,\ldots,p_n\})\sub\HH,
\end{equation}
where the embedding into $\HH$ is given by the Laurent expansions with respect to the formal parameters $(t_\bullet)$.
This morphism is $\fG$-equivariant, where $\fG=\prod_{i=1}^n\fG_i$, and $\fG_i$ is the group of changes of $t_i$
of the form $t_i\mapsto t_i+c_1t_i^2+c_2t_i^3+\ldots$.

\begin{thm}\label{big-moduli-thm} Assume that either $n\ge g\ge 1$, $n\ge 2$ and the base is $\Spec(\Z[1/2])$, or
$n=g=1$ and the base is $\Spec(\Z[1/6])$,
or $g=0$, $n\ge 2$ and the base is $\Spec(\Z)$.

\noindent
(i) The stack $\wt{\UU}_{g,n}^{ns}$ is a scheme, the morphism $\pi:\wt{\UU}_{g,n}^{ns}\to G(n-g,n)$ is 
affine of finite type and $\G_m^n$-equivariant.

\noindent
(ii) The diagonal subgroup $\G_m\sub \G_m^n$ acts on the ring of functions on each open affine $\wt{\UU}_{g,n}(S)$,
where $S\sub\{1,\ldots,n\}$, $|S|=g$, with non-negative weights.
There is a $\G_m^n$-equivariant section 
$$\si:G(n-g,n)\to \wt{\UU}_{g,n}^{ns}$$
of the morphism $\pi$, such that the locus of $\G_m$-invariant points in 
$\wt{\UU}_{g,n}^{ns}$ coincides with the image of $\si$.
These $\G_n$-invariant points correspond to the special curves
considered in Proposition \ref{special-curves-prop}.

\noindent
(iii) Now let us work over $\Q$.
Then the Krichever map \eqref{Krich-map} induces an isomorphism
$$\ov{\Kr}:\wt{\UU}^{ns}_{g,n}\simeq ASG^{ns}/\fG,$$
where the action of $\fG$ on $ASG^{ns}$ is free.
\end{thm}

\Pf . (i) The case $n=g=1$ is well known (see e.g., \cite[Thm.\ 1.2.4]{P-ainf} or \cite[Thm.\ A]{LP}), 
while the case $n\ge 3$, $g=0$ is
\cite[Thm.\ 5.1.4]{P-ainf}. The case $n=2$, $g=0$ is elementary (and can be worked out as in \cite[Lem.\ 5.1.2]{P-ainf}):
one has $\wt{\UU}^{ns}_{0,2}\simeq \A^1$, with the universal affine curve $C\setminus\{p_1,p_2\}$ given by the
equation $x_1x_2=t$ (where $t$ is the coordinate on $\A^1$).
Thus, we will assume that $n\ge g\ge 1$ and $n\ge 2$.

Since $G(n-g,n)$ is covered by affine open cells $U_S$, it is enough to prove that each
$\wt{\UU}_{g,n}(S)$ is an affine scheme of finite type over the base. As in \cite[Thm.\ 1.2.4]{P-ainf} and
\cite[Thm.\ A]{P-krich}, the main point is to find
a canonical basis of the algebra $H^0(C\setminus\{p_1,\ldots,p_n\},\OO)$ for a family of curves in
$\wt{\UU}_{g,n}(S)$ with any affine base $\Spec(R)$.

Set $D_S=\sum_{i\in S}p_i$. We start by constructing an $R$-basis of the algebra
$H^0(C\setminus D_S,\OO)$. Similarly to \cite[Thm.\ 1.2.4]{P-ainf}, using the condition $H^1(C,\OO(D_S))=0$
we choose elements 
$$f_i\in H^0(p_i+D_S), \ \ h_i\in H^0(2p_i+D_S), \text{ for } i\in S,$$ 
with $f_i\equiv \frac{1}{t_i^2}\mod \OO(D_S), h_i\equiv \frac{1}{t_i^3}\mod \OO(p_i+D_S)$,
where $t_i$ are some formal parameters at $p_i$, compatible with the chosen tangent vectors. 
The elements $(f_i,h_i)$ are defined uniquely up to the following transformations
\begin{equation}\label{f-i-h-i-transform-eq}
(h_i,f_i)\mapsto (h_i+a_if_i+b_i,f_i+c_i),
\end{equation}
for some $a_i,b_i,c_i\in R$.
The assumption that $H^1(C,\OO(D_S))=0$ implies that $H^0(C,\OO(D_S))=R$, hence,
the monomials $(f_i^m, f_i^mh_i)_{i\in S, m\ge 0}$ form an $R$-basis of
$H^0(C\setminus D_S,\OO)$ (cf. \cite[Lem.\ 1.2.1(ii)]{P-ainf}). 
By considering the polar parts, as in \cite[Lem.\ 1.2.1]{P-ainf}, 
we see that the generators $(f_i,h_i)$ should satisfy relations of the following form:
$$f_if_j=\a_{ji}h_i+\a_{ij}h_j+\ga_{ji}f_i+\ga_{ij}f_j+\sum_{k\neq i,j}c_{ij}^k f_k+a_{ij},$$
$$f_ih_j=d_{ij}f_j^2+t_{ji}h_i+v_{ij}h_j+r_{ji}f_i+\de_{ij}f_j+\sum_{k\neq i,j} e_{ij}^kf_k+b_{ij},$$
$$h_ih_j=\b_{ji}f_i^2+\b_{ij}f_j^2+\vareps_{ji}h_i+\vareps_{ij}h_j+\psi_{ji}f_i+\psi_{ij}f_j+
\sum_{k\neq i,j}l_{ij}^kf_k+u_{ij},$$
$$h_i^2=f_i^3+q_ih_if_i+r_if_i^2+u_ih_i+\sum_{j\neq i}g_i^jh_j+\pi_if_i+\sum_{j\neq i}k_i^jf_j+s_i,$$
where $i\neq j$ (the coefficients are some elements of $R$).
Since we assume that $2$ is invertible,
choosing $a_i$ and $b_i$ in \eqref{f-i-h-i-transform-eq} appropriately we can make the coefficients
$q_i$ and $r_i$ in the last equation to be zero. This fixes the ambiguity in a choice of $h_i$.
Assume now that $g\ge 2$.
Then to fix the ambiguity in a choice of $f_i$ we observe that 
by making appropriate changes $f_i\mapsto f_i+c_i$ we can make $\ga_{ii_0}=0$ for $i\neq i_0$,
$\ga_{i_0i_1}=0$ for fixed $i_0,i_1\in S$ ($i_0\neq i_1$). In the case $g=1$ we will leave the ambiguity in choosing $f_i$
for now and will fix it later.

Next, for each $j\not\in S$ we have $h^0(p_j+D_S)=2$, so we can choose $h_{S,j}\in H^0(C,\OO(p_j+D_S))$ with
the polar part $1/t_j$ at $p_j$, uniquely up to an additive constant. Let us set $D:=\sum_{i=1}^n p_i$. Then
\begin{equation}\label{big-moduli-basis-eq}
(f_i^m, f_i^mh_i, h_{S,j}^{m+1}), \ \ i\in S, j\not\in S, m\ge 0
\end{equation}
is an $R$-basis of $\OO(C\setminus D)$. Indeed, $(h_{S,j}^m)_{m\ge 1,j\not\in S}$ generate arbitrary polar parts at points $p_j$, 
$j\not\in S$, while $(f_i^m,f_i^mh_i)_{i\in S,m\ge 0}$ form a basis of $H^0(C\setminus D_S,\OO)$.

Let us define $a_{ij}(S)\in R$, where $i\in S$, $j\not\in S$, by
$$h_{S,j}\equiv \frac{a_{ij}(S)}{t_i} \mod\OO$$
at $p_i$. Then using the basis \eqref{big-moduli-basis-eq} we see that in addition to the relations satisfied by $(f_i,h_i)_{i\in S}$
we should have relations of the following form:
$$h_{S,j}h_{S,j'}=c_{j'j}(S)h_{S,j}+c_{jj'}(S)h_{S,j'}+\sum_{i\in S}a_{ij}(S)a_{ij'}(S)f_i+ const$$
$$f_ih_{S,j}=b_{ij}(S)h_{S,j}+a_{ij}(S)h_i+\sum_{l\in S} d_{ij}^l(S)f_l+ const$$
$$h_ih_{S,j}=e_{ij}(S)h_{S,j}+a_{ij}(S)f_i^2+r_{ij}(S)h_i+\sum_{l\in S} s_{ij}^l(S)f_l+ const$$
where $i\in S$, $j,j'\not\in S$.
Note that $c_{jj'}(S)=h_{S,j}(p_{j'})$, $b_{ij}(S)=f_i(p_j)$, $e_{ij}(S)=h_i(p_j)$.
Using these relations we can get rid of the ambiguity in adding a constant to each $h_{S,j}$
by requiring that $d_{i_0j}^{i_0}(S)=0$ for fixed $i_0\in S$.
Also, in the case $g=1$ we can now fix the ambiguity in adding a constant to $f_i$ by requiring $b_{ij_0}(S)=0$ for
a fixed $j_0\not\in S$ (which exists since $|S|=1$ and $n\ge 2$). 

As usual, the Buchberger's algorithm gives a system of equations on the constants in the relations between the generators
$f_i,h_i,h_{S,j}$, which is equivalent to \eqref{big-moduli-basis-eq} being a basis. Thus, we get a morphism from
$\wt{\UU}_{g,n}(S)$ to an affine scheme $S_{GB}$ of finite type over $\Z[1/6]$. The remainder of the proof
is similar to that of \cite[Thm.\ 1.2.4]{P-ainf}: 
starting from Groebner relations of the above form we construct a family of
curves with required properties, which gives an inverse morphism $S_{GB}\to \wt{\UU}_{g,n}(S)$.

\noindent
(ii) The action of $\la\in \G_m^n$ on the coordinates on $\wt{\UU}_{g,n}(S)$ is induced by the rescalings
$$f_i\mapsto \la_i^2f_i, \ \ h_i\mapsto \la_i^3 h_i, \ \ h_{S,j}\mapsto \la_jh_{S,j}.$$
From this we see that the diagonal action of $\G_m$ doesn't change $(a_{ij})$ and acts with positive weights on all the other 
coordinates. In particular, the $\G_m$-invariant points correspond to the $(C,p_\bullet,v_\bullet)$ such that 
$C\setminus \{p_1,\ldots,p_n\}$ is given by the equations \eqref{Gm-inv-curves-eq}.
So these are exactly the equations of the special curves 
considered in Proposition \ref{special-curves-prop}. 

Conversely, from part (i), we see that $\wt{\UU}_{g,n}(S)$ can be identified
with the affine scheme parametrizing commutative algebras with generators $(f_i, h_i, h_{S,j})$ and relations of the form
specified in (i), such that the elements  \eqref{big-moduli-basis-eq} form a basis. Thus, by Proposition
\ref{special-curves-prop}, the family of special curves over $U_S$ gives a morphism
$$\si:U_S\to \wt{\UU}_{g,n}(S).$$ 
It remains to check that this morphism is a section of the projection $\pi:\wt{\UU}_{g,n}(S)\to U_S$.
By definition, we have to prove that for $C=C(a_{ij})$ as in Proposition \ref{special-curves-prop}, where $(a_{ij})$
is the matrix with $i\in S$, $j\not\in S$, the kernel of the map 
$$H^0(C,\OO_C(p_1+\ldots+p_n)/\OO_C)\to H^1(C,\OO_C)$$
gets identified with the subspace $W\sub k^n$ spanned by $(\be_j+\sum_{i\in S} a_{ij}\be_i)_{j\not\in S}$,
where $H^0(C,\OO(p_i)/\OO)$ is trivialized using the rational function $x_i$ that has a pole of order $1$ at $p_i$.
But this follows from the fact that $x_j+\sum_{i\in S} a_{ij}x_i=h_{S,j}$ defines a regular section in
$H^0(C,\OO_C(p_1+\ldots+p_n))$.

\noindent
(iii) We have seen that $\wt{\UU}^{ns}_{g,n}$ is the union of $\wt{\UU}^{ns}_{g,n}(\ba)$, over $\ba$ consisting of $0$'s and
$1$'s. Similarly, $ASG^{ns}$ is the union of $ASG^\ba$ (where $ASG^\ba=ASG\cap SG^\ba$ with
$SG^\ba$ being the open cell of the Sato Grassmannian defined in \cite[Sec.\ 1.3]{P-krich}), 
over the same set of $\ba$, and the Krichever map is
compatible with these open coverings. Hence, the assertion follows from \cite[Thm.\ B]{P-krich}.
\ed

Note that Theorem A is a part of Theorem \ref{big-moduli-thm}.

\begin{rems}
1. If $C$ has arithmetic genus $0$ then the vanishing of $H^1(\OO(p_1+\ldots p_n))$ is automatic, so the moduli scheme
$\wt{\UU}_{0,n}^{ns}$ is exactly the space $\wt{\UU}_{0,n}[\psi]$ of $\psi$-prestable curves (with tangent vectors
at the marked points) considered in \cite[Sec.\ 5.1]{P-ainf}. This case of Theorem \ref{big-moduli-thm}(i)(ii) was considered in \cite[Thm.\ 5.1.4]{P-ainf}. Note that one of the GIT quotients of $\wt{\UU}_{0,n}^{ns}$
by $\G_m^n$ is the space of Boggi-stable (or $\psi$-stable) curves studied in
\cite{Boggi}, \cite[Sec.\ 4.2.1]{FS} and \cite[Sec.\ 7.2]{GJM}.

\noindent
2. The schemes $\wt{\UU}_{g,n}^{ns}$ can be reducible. For example, by \cite[Thm.\ (11.10)]{Pin},
if $C$ is the union of $n$ generic lines through one point in $\P^{n-3}$, then $C$ is not smoothable for $n\ge 15$.
It is easy to see that equipping each component of $C$ with a marked point we get a curve $(C,p_1,\ldots,p_n)$ satisfying
$H^1(C,\OO(p_1+\ldots,p_n))=0$. Thus, we deduce that $\wt{\UU}_{3,n}^{ns}$, for $n\ge 15$, has a component with
nonsmoothable curves.
\end{rems}

\begin{defi}\label{U-prime-def}
Let us denote by $\UU^{ns, \prime}_{g,n}$ (resp., $\wt{\UU}^{ns, \prime}_{g,n}$)
the moduli stack, defined exactly like $\UU_{g,n}^{ns}$ (resp., $\wt{\UU}^{ns}_{g,n}$), but without the condition
of ampleness of $\OO(p_1+\ldots+p_n)$.
\end{defi}

\begin{prop}\label{U-prime-prop} 
There is a natural morphism $\UU^{ns, \prime}_{g,n}\to \UU^{ns}_{g,n}$ sending a curve
$(C,p_1,\ldots,p_n)$ to the curve $(\ov{C},\ov{p}_1,\ldots,\ov{p}_n)$, where $\ov{C}$ is 
the image of $C$ under the morphism to a projective space
induced by the linear system $|\OO_C(N(p_1+\ldots+p_n))|$ for $N\gg 0$.
\end{prop}

\Pf . We will construct a $\G_m^n$-equivariant morphism
$$\wt{\UU}^{ns, \prime}_{g,n}\to\wt{\UU}_{g,n}^{ns},$$
as sketched in \cite[Rem.\ 1.7.2.2]{P-krich}.
Let us restrict to the open substack where $h^1(\OO(\sum_{i\in S} p_i))=0$, for a fixed subset
$S\sub\{1,\ldots,n\}$, $|S|=g$ (it will be clear that our morphisms are the same on the intersections).

Let $R$ be a commutative ring, and let $(C,p_\bullet,v_\bullet)$ be a family in $\wt{\MM}_{1,n}(m)(R)$.
Consider the $R$-algebra 
$$A=H^0(C\setminus \{p_1,\ldots,p_n\},\OO),$$ 
equipped with the increasing filtration $F_m=H^0(C,\OO(mD))$, where $D=p_1+\ldots+p_n$.
Now set 
$$\ov{C}:=\Proj(\RR(A)),$$ 
where $\RR(A)=\bigoplus_m F_m$ is the corresponding Rees algebra. 
As in the proof of Theorem \ref{big-moduli-thm}, we construct a canonical basis of the algebra $A$
(here we use the assumption that $h^1(\OO(\sum_{i\in S} p_i))=0$).
Then, following the argument of the same proof, we
use the relations between the canonical 
generators of $A$ to show that $\ov{C}$ is equipped with marked points and tangent vectors and defines 
a family in $\wt{\UU}_{1,n}$.
\ed

\subsection{Gluing morphism}\label{gluing-sec}

Let $\wt{\CC}_{g,n}^{ns}\to \wt{\UU}_{g,n}^{ns}$ denote the universal {\it affine} curve, i.e.,
the complement to the sections $p_1,\ldots,p_n$ in the universal curve.
By definition, the stack $\wt{\CC}_{g,n}^{ns}$ classifies the data $(C,p_1,\ldots,p_n,v_1,\ldots,v_n;q)$,
where $(C,p_\bullet,v_\bullet)$ is in $\wt{\UU}_{g,n}^{ns}$ and $q$ is a point of $C$, different from the marked points
$p_1,\ldots,p_n$ (where $C$ can be singular at $q$).
In the case when the marked points are in bijection with a finite set $I$ we will use the notation
$\wt{\UU}_{g,I}^{ns}$ and $\wt{\CC}_{g,I}^{ns}$ for these moduli stacks.

\begin{ex}
In the case $n=g=1$ we need to invert $6$ to ensure that the stack $\wt{\UU}^{ns}_{1,1}$ is a scheme.
However, it is easy to see that already the stack $\wt{\CC}_{1,1}^{ns}\times\Spec(\Z[1/2])$ is a scheme.
Indeed, starting with a family $(C,p,v;q)$ in $\wt{\CC}_{1,1}^{ns}(R)$, where $R$ is a commutative ring,
we can normalize functions $f\in H^0(C,\OO(2p))$, $h\in H^0(C,\OO(3p))$, with the Laurent
expansions $f=\frac{1}{t^2}+\ldots$, $g=\frac{1}{t^3}+\ldots$ at $p$ (where the local parameter $t$ is compatible with the vector
field $v$) by the conditions $f(q)=h(q)=0$, so that the only remaining ambiguity is $h\mapsto h+cf$.
We can fix this ambiguity by requiring that 
$$h^2-f^3=\a f^2+\b h+\ga f$$
for uniquely defined constants $\a,\b,\ga\in R$ (note that there is no $fh$ term in the right-hand side). 
This gives an isomorphism $\wt{\CC}_{1,1}^{ns}\simeq \A^3$ over $\Z[1/2]$.
\end{ex}

\begin{prop} Assume $n\ge g\ge 1$, $n\ge 2$, and let us work over $\Z[1/2]$ (in the case $g=0$ we can work over $\Z$).

\noindent
(i) For every partition $[1,n]=I\sqcup J$ into nonempty subsets
and a pair of numbers $h\le |I|, k\le |J|$, such that $h+k=g$, there is a natural {\it gluing morphism}
\begin{equation}
\rho^{h,k}_{I,J}:\wt{\CC}_{h,I}^{ns}\times \wt{\CC}_{k,J}^{ns}\to \wt{\UU}_{g,n}^{ns},
\end{equation}
sending a pair of curves $C_I$, $C_J$ (equipped with the marked points, tangent vectors and
with the extra points $q_I$, $q_J$) to the curve $C=C_I\cup C_J$,
where $C_I$ and $C_J$ are glued transversally, so that the points $q_I$ and $q_J$ are identified.
The curve $C$ is equipped with $n$ marked points (and tangent vectors at them), 
so that the points indexed by $I$ come
from the marked points on $C_I$, while those indexed by $J$ come from $C_J$.

\noindent
(ii) The morphism $\rho^{h,k}_{I,J}$ is compatible with the projections to the Grassmannians and with the morphism
$$G(|I|-h,|I|)\times G(|J|-k,|J|)\to G(n-g,n)$$
sending the pair of subspaces $(V,V')$ to $V\oplus V'$.

\noindent
(iii) The morphism $\rho^{h,k}_{I,J}$ is a closed embedding that factors through the closed subscheme 
$Z_{I,J}\sub\wt{\UU}_{g,n}^{ns}$, given by the conditions
\begin{enumerate}
\item
$H^1\bigl(C,\OO(\sum_{i\in S} p_i)\bigr)\neq 0$ for all $S\sub [1,n]$ such that $|S|=g$ and either $|S\cap I|<h$ or $|S\cap J|<k$;
\item
for every $S\sub [1,n]$ and $s,t\in [1,n]$, such that $|S|=g$, $|S\cap I|=h$, $|S\cap J|=k$, and either $s\in S\cap I$, $t\in S\cap J$
or $s\in S\cap J$, $t\in S\cap I$, the morphism 
$$H^0\bigl(C,\OO(2p_s+\sum_{i=1}^n p_i)\bigr)\to H^0(C,\OO(p_t)/\OO)\oplus\bigoplus_{i\not\in S}H^0(C,\OO(p_i)/\OO)$$ 
has rank $\le n-g$.
\end{enumerate}
Note that the nonvanishing of $H^1$ in (1) can also be expressed as a degeneracy locus of the morphism of vector bundles
over $\wt{\UU}_{g,n}^{ns}$, so we have a natural subscheme structure on $Z_{I,J}$.
Furthermore, there exists a retraction morphism from $Z_{I,J}$ onto the image of $\rho^{h,k}_{I,J}$.
\end{prop} 

\Pf . (i,ii) The fact that $C_I$ and $C_J$ are glued transversally means that there is an exact sequence
$$0\to \OO_C\to \OO_{C_I}\oplus\OO_{C_J}\to \OO_q\to 0.$$
This leads to an isomorphism 
\begin{equation}\label{glued-H1-decomposition}
H^1(\OO_C)\simeq H^1(\OO_{C_I})\oplus H^1(\OO_{C_J}),
\end{equation}
so the arithmetic genus of $C$ is $h+k=g$. Similarly, the exact sequence
$$0\to \OO_C(p_1+\ldots+p_n)\to \OO_{C_I}(\sum_{i\in I}p_i)\oplus\OO_{C_J}(\sum_{j\in J}p_j)\to \OO_q\to 0$$
shows that
$$H^1\bigl(\OO_C(p_1+\ldots+p_n)\bigr)=H^1\bigl(\OO_{C_I}(\sum_{i\in I}p_i)\bigr)\oplus 
H^1\bigl(\OO_{C_J}(\sum_{j\in J}p_j)\bigr)=0$$
(since the constants in $H^0\bigl(\OO_{C_I}(\sum_{i\in I}p_i)\bigr)$ surject onto $H^0(\OO_q)$).
Similar argument works in families, so our morphism is well-defined.
The compatibility with the morphism of Grassmannians follows from the decomposition \eqref{glued-H1-decomposition},
which is compatible with the similar decomposition of $H^1(\OO_C(p_1+\ldots+p_n)/\OO_C)$.

\noindent
(iii) For any subset $S\sub [1,n]$ we have an exact sequence
$$0\to \OO_C(\sum_{i\in S}p_i)\to \OO_{C_I}(\sum_{i\in S\cap I}p_i)\oplus\OO_{C_J}(\sum_{j\in S\cap J}p_j)\to \OO_q\to 0$$
which gives an isomorphism
$$H^1\bigl(C,\OO_C(\sum_{i\in S}p_i)\bigr)\simeq H^1\bigl(C_I,\OO(\sum_{i\in S\cap I}p_i)\bigr)\oplus 
H^1\bigl(C_J,\OO(\sum_{j\in S\cap J}p_j)\bigr).$$
Thus, if $|S\cap I|<h$ then $H^1\bigl(C_I,\OO(\sum_{i\in S\cap I}p_i)\bigr)\neq 0$, so that 
$H^1\bigl(C,\OO_C(\sum_{i\in S}p_i)\bigr)\neq 0$.
Similarly, we get the nonvanishing of $H^1\bigl(C,\OO_C(\sum_{i\in S}p_i)\bigr)$ if $|S\cap J|<k$.
Assume now that $S$ is as in condition (2), $s\in S\cap I$ and $t\in S\cap J$. 
Using the exact sequence 
$$0\to \OO_C(2p_s+\sum_{i=1}^n p_i)\to \OO_{C_I}(2p_s+\sum_{i\in I}p_i)\oplus\OO_{C_J}(\sum_{j\in J}p_j)\to \OO_q\to 0$$
we see that the degeneracy required in (2) is equivalent to the condition that the morphism
$$H^0\bigl(C_I,\OO(2p_s+\sum_{i\in I}p_i)\bigr)\oplus H^0\bigl(C_J,\OO(\sum_{j\in J}p_j)\bigr)\to 
H^0(\OO_q)\oplus H^0(\OO(p_t)/\OO)\oplus\bigoplus_{i\not\in S}H^0(\OO(p_i)/\OO)$$
has rank $\le n-g+1$. But this follows from the fact that the composition of this morphism with the projection to 
$H^0(\OO(p_t)/\OO)\oplus\bigoplus_{j\in J\setminus S}H^0(\OO(p_j)/\OO)$
has rank $\le |J|-k$. Indeed, this composition factors through a map
$$H^0(C_J,\OO(\sum_{j\in J}p_j))\to H^0(\OO(p_t)/\OO)\oplus\bigoplus_{j\in J\setminus S}H^0(\OO(p_j)/\OO)$$
whose cokernel is $H^1\bigl(C_J,\OO(\sum_{j\in S\cap J\setminus \{t\}})\bigr)\neq 0$.
This shows that our morphism factors through the subscheme $Z_{I,J}$.

Next, we are going to construct a retraction 
\begin{equation}\label{r-retraction-map}
r:Z_{I,J}\to \wt{\CC}_{h,I}^{ns}\times \wt{\CC}_{k,J}^{ns}.
\end{equation}
It is enough to construct compatible morphisms on all the affine opens $Z_{I,J}\cap \wt{\UU}_{g,n}(S)$,
where $S\sub [1,n]$, $|S|=g$. By condition (1), the intersection is nonempty only when $|S\cap I|=h$ and $|S\cap J|=k$.
Let $(C,p_\bullet,v_\bullet)$ be the restriction of the universal family to $Z_{I,J}\cap \wt{\UU}_{g,n}(S)$.

Recall that over $\wt{\UU}_{g,n}(S)$ we have generators $f_i$, $h_i$, $h_{S,j}$ of $\OO(C\setminus D)$, where 
$D=p_1+\ldots+p_n$ (see the proof of Theorem \ref{big-moduli-thm}). 
Let us consider their restrictions over $Z_{I,J}$ (denoted in the same way).
We claim that for every $i\in S\cap I$ and $i'\in I\setminus S$
the functions $f_i$, $h_i$ and $h_{S,i'}$ are regular at any $p_j$, where $j\in S\cap J$.
Indeed, for $h_{S,i'}$ this follows from the exact
sequence
$$H^0\bigl(C,\OO(p_{i'}+\sum_{i\in S}p_i)\bigr)\to H^0(C,\OO(p_{i'})/\OO)\to 
H^1\bigl(C,\OO(p_{i'}+\sum_{i\in S\setminus\{j\}}p_i)\bigr)\to 0$$
since by condition (1) the first arrow is zero (due to the way we represent the nonvanishing of $H^1$ as 
a degeneracy locus).
Since over $\wt{\UU}_{g,n}(S)$, for $i\in S\cap I$, we have an exact sequence
$$0\to H^0\bigl(C,\OO(2p_i+\sum_{i'\in S}p_{i'})\bigr)\to H^0\bigl(C,\OO(2p_i+\sum_{i'=1}^n p_{i'})\bigr)\to \bigoplus_{i'\not\in S}H^0(\OO(p_{i'})/\OO)
\to 0,$$
condition (2) implies the vanishing of the morphism
$$H^0\bigl(C,\OO(2p_i+\sum_{i'\in S}p_{i'})\bigr)\to H^0(\OO(p_j)/\OO)$$
for any $j\in S\cap J$. This shows that $f_i$ and $h_i$ have no poles along $p_j$ for such $j$,
proving our claim.

Now we construct two families of curves over $Z_{I,J}$. Let us set
$$A_I:=\OO(C\setminus\{p_i\ |\ i\in I\}), \ \ A_J:=\OO(C\setminus\{p_i\ |\ i\in J\}),$$
$$C_I:=\Proj \RR(A_I), \ \ C_J:=\Proj \RR(A_J),$$
where $\RR(A_I)$ (resp., $\RR(A_J)$) are the Rees algebras associated with the filtrations by the order of pole
along $\sum_{i\in I} p_i$ (resp., $\sum_{j\in J}p_j$). Note that the algebra $A_I$ is a free module over
$\OO(Z_{I,J})$ with the basis 
$$(f_i^m, f_i^mh_i, h_{S,i'}^{m+1}), \ \ i\in S\cap I, i'\in I\setminus S, m\ge 0.$$
Indeed, this
is checked easily by considering polar parts at $p_i$ with $i\in I$, similarly to checking that \eqref{big-moduli-basis-eq} 
is a basis of $\OO(C\setminus D)$. This implies that the generators $(f_i,h_i,h_{S,i'})$, where $i\in S\cap I$,
$i'\in I\setminus S$, satsify similar relations as the full set of generators of $\OO(C\setminus D)$, and so as in the proof of
\cite[Thm.\ 1.2.4]{P-ainf}, 
we get that $C_I$ is a curve of arithmetic genus $|S\cap I|=h$, that extends to a family in $\wt{\UU}_{h,I}^{ns}$.
Similar argument works for $C_J$, so we get a morphism
$$\ov{r}:Z_{I,J}\to \wt{\UU}_{h,I}^{ns}\times \wt{\UU}_{k,J}^{ns}.$$
Note that we have the natural morphisms
$C\to C_I$ and $C\to C_J$. Let us consider the compositions
$$q:Z_{I,J}\rTo{p_{j_0}} C\to C_I, \ \ q':Z_{I,J}\rTo{p_{i_0}} C\to C_J,$$
for some fixed indices $i_0\in I$, $j_0\in J$.
These allow to lift the morphism $\ov{r}$ to the required morphism \eqref{r-retraction-map}.
One can immediately check that $r\circ \rho_{I,J}^{h,k}$ is the identity morphism.
This implies that $\rho_{I,J}^{h,k}$ is a closed embedding.
\ed

\begin{ex}
In the case $g=0$ the subscheme $Z_{I,J}$ coincides with the whole space $\wt{\UU}_{0,n}^{ns}$. Thus,
in this case the embedding $\rho^{0,0}_{I,J}$ admits a retraction defined on $\wt{\UU}_{0,n}^{ns}$.
\end{ex}

\begin{rem} One can also consider the gluing morphism that associates with a curve $(C,p_1,\ldots,p_n)$ of arithmetic genus $g-1$, equipped with extra two points $q_1\neq q_2$ (possibly singular), different from $p_1,\ldots,p_n$, the curve 
$(\ov{C},p_1,\ldots,p_n)$ of arithmetic genus $g$,
where $\ov{C}$ is obtained by transversally identifying $q_1$ and $q_2$ on $C$. In order to guarantee that
$H^1(\ov{C},\OO(p_1+\ldots+p_n))=0$ one has to assume that $H^1(C,\OO(p_1+\ldots+p_n))=0$ and in addition
the morphism 
$$\ev_{q_2}-\ev_{q_1}:H^0(C,\OO(p_1+\ldots+p_n))\to k,$$ is surjective (where $\ev_{q_i}$ is the evaluation at $q_i$).
\end{rem}

\subsection{Case $g=1$, $n=2$}\label{g=1-n=2-sec}
In this section we work over $\Spec(\Z[1/6])$.
The scheme $\wt{\UU}_{1,2}^{ns}$ is glued from two affine open pieces $\UU_1$ and $\UU_2$
determined by the conditions $H^1(C,\OO(p_1))=0$ and $H^1(C,\OO(p_2))=0$, respectively.

Note that we have $\UU_1=\wt{\UU}^{ns}_{1,2}(1,0)$. The latter moduli space was described explicitly
in \cite[Sec.\ 3.1]{P-krich} as the affine $4$-space with coordinates $(a,b,e,\pi)$. Let us
rename these coordinates on $\UU_1$ as $a_{12}, b_{12}, e_{12}, \pi_1$. Thus, the universal affine curve $C\setminus\{p_1,p_2\}$ over $\UU_1$ is given by the equations 
\begin{equation}\label{g-1-n-2-univ-curve-U1-eq1}
h_1^2=f_1^3+\pi_1 f_1+s_1,
\end{equation}
\begin{equation}\label{g-1-n-2-univ-curve-U1-eq2}
f_1h_{12}=a_{12}h_1+b_{12}h_{12}+a_{12}e_{12},
\end{equation}
\begin{equation}\label{g-1-n-2-univ-curve-U1-eq3}
h_1h_{12}=a_{12}f_1^2+e_{12}h_{12}+a_{12}b_{12}f_1+a_{12}(\pi_1+b_{12}^2),
\end{equation}
where $s_1=e_{12}^2-b_{12}(\pi_1+b_{12}^2)$. 
Note that the projection $\UU_1\to \A^1\sub\P^1$ is given by the coordinate $a_{12}$.

Let us denote by $f_2,h_2, h_{21}$ the generators of the algebra of functions on the universal
affine curve $C\setminus\{p_1,p_2\}$ over $\UU_2$, so
that $h_{21}\in H^0(C,\OO(p_1+p_2))$ and $h_{21}\equiv 1/t_1$ at $p_1$.
Let also $a_{21}, b_{21}, e_{21}, \pi_2$ be the coordinates on $\UU_2$ similar to those on $\UU_1$.
Note that over $\UU_1\cap \UU_2$
the function $h_{21}-\frac{1}{a_{12}}h_{12}$ is constant along the fibers of the projection to the base.
Hence, $h_{21}\equiv \frac{1}{a_{12} t_2}$ at $p_2$, so we have
$$a_{21}=\frac{1}{a_{12}}.$$

\begin{lem}\label{U-1-cap-U-2-lem} 
Over $\UU_1\cap \UU_2$ one has
$$h_{21}=a_{21}h_{12}, \ f_2=h_{12}^2-a_{12}^2f_1-a_{12}^2b_{12}, \ 
h_2=h_{12}^3-a_{12}^3h_1-3a_{12}^2b_{12}h_{12}-2a_{12}^3e_{12},$$
$$b_{21}=a_{12}^2 b_{12}, \ \ e_{21}=a_{12}^3e_{12}, \ \ \pi_2=a_{12}^4\pi_1, \ \ s_2=a_{12}^6s_1.$$ 
\end{lem}

\Pf . Note that $\UU_1\cap \UU_2$ is the open subset in $\UU_1$ where $a_{12}$ does not vanish.
On this open subset we can express $h_1$ in terms of $f_1$ and $f_{12}$ using \eqref{g-1-n-2-univ-curve-U1-eq2}:
\begin{equation}\label{h_1-U-1-2-eq}
h_1=\frac{f_1h_{12}}{a_{12}}-\frac{b_{12}h_{12}}{a_{12}}-e_{12}.
\end{equation}
Substituting into \eqref{g-1-n-2-univ-curve-U1-eq3} we get
$$f_1^2-(\frac{h_{12}^2}{a_{12}^2}-b_{12})f_1+\frac{b_{12}h_{12}^2}{a_{12}^2}+
2e_{12}\frac{h_{12}}{a_{12}}+\pi_1+b_{12}^2=0$$
which is the single equation defining the universal $C\setminus\{p_1,p_2\}$ 
(\eqref{g-1-n-2-univ-curve-U1-eq1} follows from this and \eqref{h_1-U-1-2-eq}).
Note that this is a quadratic equation in $f_1$.
Therefore, we can define an involution of $C\setminus\{p_1,p_2\}$ over $\UU_1\cap \UU_2$ by 
$$(h_{12},f_1)\mapsto (h_{12},\frac{h_{12}^2}{a_{12}^2}-b_{12}-f_1).$$
We claim that this involution extends to an involution $\tau$ of $C$ permuting $p_1$ and $p_2$.
Indeed, this immediately follows from the fact that it preserves the filtration by the degree of pole along $p_1+p_2$
(so that $\deg(h_{12})=1$, $\deg(f_1)=2$), together with the observation that $\tau^*(f_1)$ has a pole at $p_2$. 

Note that 
\begin{equation}\label{tau*-f-1-eq}
\tau^*(f_1)=\frac{h_{12}^2}{a_{12}^2}-b_{12}-f_1,
\end{equation}
which has the expansion $\frac{1}{a_{12}^2 t_2^2}+\ldots$ at $p_2$. 
Similarly 
$$\tau^*(h_1)=\frac{h_{12}}{a_{12}}\tau^*(f_1)-\frac{b_{12}h_{12}}{a_{12}}-e_{12}\equiv \frac{1}{a_{12}^3t_2^3}+\ldots$$
at $p_2$. Now the equation 
$$\tau^*(h_1)^2=\tau^*(f_1)^3+\pi_1\tau^*(f_1)+s_1$$
(obtained from \eqref{g-1-n-2-univ-curve-U1-eq1}), together with the definition of $(f_2,h_2)$, shows that
\begin{equation}\label{tau-f-h-eq}
\tau^*(f_1)=\frac{1}{a_{12}^2}f_2, \ \ \tau^*(h_1)=\frac{1}{a_{12}^3}h_2,
\end{equation}
$$\pi_2=a_{12}^4\pi_1, \ \ s_2=a_{12}^6s_1.$$
We also have 
$$b_{21}=f_2(p_1)=a_{12}^2(\tau^*f_1)(p_1)=a_{12}^2f_1(\tau(p_1))=a_{12}^2f_1(p_2)=a_{12}^2b_{12}.$$
Similarly, we get that $e_{21}=a_{12}^3e_{21}$. 
Finally, we get the required formulas for $f_2$ and $h_2$ by using \eqref{tau*-f-1-eq},
\eqref{tau-f-h-eq} and \eqref{h_1-U-1-2-eq}.
\ed

The above lemma shows that the functions $b_{12},e_{12},\pi_1$ (resp., $b_{21},e_{21},\pi_{21}$) defined on $\UU_1$ (resp., $\UU_2$) actually extend to regular functions on the entire moduli space.

\begin{prop}\label{U-1-2-prop} Let us work over $\Z[1/6]$.
The scheme $\wt{\UU}_{1,2}^{ns}$ is isomorphic as a $\P^1$-scheme 
to the total space of the vector bundle $\OO(-2)\oplus\OO(-3)\oplus\OO(-4)$ over $\P^1$.  
This isomorphism is $\G_m^2$-equivariant, where we use 
the natural action of $\G_m^2\sub\GL_2$ on $\P^1$ and its standard lifting to $\OO(i)$.
\end{prop}

\Pf . The morphism $\pi:\wt{\UU}_{1,2}^{ns}\to\P^1$ is given by $(1:a_{12})$ on $\UU_1$
and by $(a_{21}:1)$ on $\UU_2$, where $a_{12}=a_{21}^{-1}$ on the intersection.
Since $\wt{\UU}_{1,2}^{ns}$ is glued from the open subsets $\UU_1$ and $\UU_2$, using the identifications
of $\UU_1$ and $\UU_2$ with $\A^4$ and the transition formulas from Lemma \ref{U-1-cap-U-2-lem},
we obtain that $\wt{\UU}_{1,2}^{ns}$ is isomorphic to the 
subscheme of $\P^1\times\A^6$ given by the
equations
$$t_1^2b_{21}=t_2^2b_{12}, \ \ t_1^3e_{21}=t_2^3e_{12}, \ \ t_1^4\pi_2=t_2^4\pi_1,$$
where $(t_1:t_2)$ are homogeneous coordinates on $\P^1$.
Hence, $b_{12}/t_1^2$, $e_{12}/t_1^3$ and $\pi_1/t_1^4$ extend naturally to regular sections of
$\pi^*\OO(-2)$, $\pi^*\OO(-3)$ and $\pi^*\OO(-4)$ respectively. This gives a morphism from
$\wt{\UU}_{1,2}^{ns}$ to the total space of $\OO(-2)\oplus\OO(-3)\oplus\OO(-4)$ over $\P^1$.
Since it restricts to isomorphisms over the open subsets $t_1\neq 0$ and $t_2\neq 0$, it is an isomorphism.
\ed

\begin{rem} Under the isomorphism of Proposition \ref{U-1-2-prop}, 
the $\G_m$-invariant points in $\wt{\UU}_{1,2}^{ns}$ (see Theorem \ref{big-moduli-thm}(ii))
get identified with the zero section in the total space of $\OO(-2)\oplus\OO(-3)\oplus\OO(-4)$ over $\P^1$.
Over $\P^1\setminus\{0,\infty\}$ the corresponding curve is the tacnode, while at $0$ and $\infty$ we get the union of the
genus $1$ cuspidal curve and of the projective line, joined to the cusp transversally.
\end{rem}



In the remainder of this section we work over an algebraically closed field $k$ of characteristic $\neq 2,3$.

\begin{cor}\label{U-1-2-algebra-cor}
The graded algebra
$$A=\bigoplus_{n\ge 0}H^0(\wt{\UU}_{1,2}^{ns}, \pi^*\OO(n))=
\bigoplus_{n\ge 0}H^0(\P^1,S(\OO(2))\ot S(\OO(3))\ot S(\OO(4))(n))$$ 
can be identified with the $k[t_1,t_2]$-subalgebra of $k[t_1,t_2,x,y,z]$ (where $t_1,t_2,x,y,z$ are independent variables), 
such that the $n$th graded component
$A_n$ is spanned by the monomials
$$t_1^it_2^jx^ky^lz^m \ \ \text{with } i+j=2k+3l+4m+n.$$
This identification is compatible with the $\G_m^2$-actions, where $x,y,z$ and $\G_m^2$-invariant,
while $t_1$ and $t_2$ have $\G_m^2$-weights $(1,0)$ and $(0,1)$, respectively.
\end{cor}

\Pf . We use the $\G_m^2$-invariant regular sections of $\pi^*\OO(-2)$, $\pi^*\OO(-3)$ and $\pi^*\OO(-3)$,
$$x=b_{12}/t_1^2, \ \ y=e_{12}/t_1^3, \ \ z=\pi_1/t_1^4.$$
\ed

Using the above description we can describe the GIT quotients $\wt{\UU}^{ns}_{1,2}\sslash_\chi \G_m^2$ with
respect to the $\G_m^2$-linearizations on the line bundle $\OO(1)$ on $\wt{\UU}^{ns}_{1,2}$, which differ
from the standard $\G_m^2$-equivariant structure by (rational) characters 
$\chi(\la_0,\la_1)=\la_0^u\la_1^v$ of $\G_m^2$ (if $u$ and $v$ are fractional this means that we really work
with $\OO(N)$ for some $N$).

We have
$$\wt{\UU}^{ns}_{1,2}\sslash_\chi \G_m^2=\Proj(A(u,v)),$$
where $A(u,v)\sub A$ is the corresponding subalgebra of invariants with respect to the
$\chi$-twisted $\G_m^2$-action on the algebra $A$:
$$A(u,v):=\bigoplus_{n\ge 0}(A_n\ot\chi^{-n})^{\G_m^2}$$
(here, if $\chi$ is not integral, we have to pass to a Veronese subalgebra of $A$).
Using Corollary \ref{U-1-2-algebra-cor}
we see that the $A(u,v)_n$ is spanned by the monomials
\begin{equation}\label{xy-twisted-invariants-eq}
t_1^{nu}t_2^{nv}x^ky^lz^m \ \ \text{ with } n(u+v-1)=2k+3l+4m.
\end{equation}

Note that $A(u,v)$ reduces to constants if either $u<0$ or $v<0$ or $u+v<1$.
In the case $u\ge 0$, $v\ge 0$, $u+v=1$, the algebra $A(u,v)$ has a basis of monomials
$t_1^{nu}t_2^{nv}$, so it is isomorphic to an algebra of polynomials in one variable $t_1^{Nu}t_2^{Nv}$,
where $N>0$ is minimal such that $Nu,Nv$ are integers. Thus, in this case 
the GIT quotient reduces to a point. 

\begin{prop} Assume $u\ge 0$, $v\ge 0$, $u+v>1$. Then
$A(u,v)$ is isomorphic to a Veronese subalgebra in
$k[x,y,z]$, where $\deg(x)=2$, $\deg(y)=3$, $\deg(z)=3$. Hence, for $u+v>1$ the GIT quotient 
$\wt{\UU}^{ns}_{1,2}\sslash_\chi \G_m^2$
is isomorphic to the weighted projective plane $\P(2,3,4)$.
The $\chi$-unstable locus in $\wt{\UU}^{ns}_{1,2}$ is the union of the locus of $\G_m$-invariant points and of
$\pi^{-1}(t_1^{Nu}=0)\cup \pi^{-1}(t_2^{Nv}=0)$, where $N>0$ is such that $Nu,Nv\in\Z$.
\end{prop}

\Pf . Let $N>0$ be minimal such that $Nu, Nv$ are integers.
We associate with each monomial \eqref{xy-twisted-invariants-eq} the corresponding monomial
$x^ky^lz^m$. This gives an isomorphism of $A(u,v)$ with the subalgebra of $k[x,y,z]$ spanned by all such monomials
such that $n:=\frac{2k+3l+4m}{u+v-1}$ belongs to $N\Z$. Let $n_0$ be minimal such that $n>0$. Then
this subalgebra is precisely the Veronese subalgebra corresponding to $n_0(u+v-1)$. The identification of the unstable
locus follows from the form of the monomials \eqref{xy-twisted-invariants-eq}, since the locus of $\G_m$-invariant points 
coincides with the locus where all the sections $x,y,z$ vanish.
\ed



\subsection{Case $g=1$, aribitrary $n\ge 2$: connection to Smyth's moduli}\label{g=1-any-n-sec}

As before, we work over $\Z[1/6]$.

Recall that for integers $1\le m<n$ Smyth defined in \cite{Smyth-I} the moduli stack $\ov{\MM}_{1,n}(m)$ of $n$-pointed
$m$-stable curves of arithmetic genus $1$, parametrizing curves $(C,p_1,\ldots,p_n)$ of arithmetic genus $1$
with $n$ distinct smooth marked points such that 
\begin{itemize}
\item $C$ has only nodes and elliptic $l$-fold points, with $l\le m$, as singularities;
\item if $E\sub C$ is any connected subcurve of arithmetic genus $1$ then 
$|E\cap \ov{C\setminus E}|+|E\cap\{p_1,\ldots,p_n\}|>m$;
\item $H^0(C,\TT_C(-p_1-\ldots-p_n))=0$, where $\TT_C$ is the tangent sheaf.
\end{itemize}
Smyth showed that $\ov{\MM}_{1,n}(m)$ is a proper irreducible Deligne-Mumford stack. 

\begin{prop}\label{Smyth-moduli-prop}
Assume $m\ge \frac{n-1}{2}$. Then there exists a morphism
$$\ov{\MM}_{1,n}(m)\to \UU_{1,n}^{ns}$$
extending the obvious map on the locus of smooth curves.
\end{prop}

\Pf . 
It is enough to check that $\ov{\MM}_{1,n}(m)$ is an open substack of $\UU^{ns, \prime}_{1,n}$ (see 
Definition \ref{U-prime-def})
for these values of $m$. Indeed, then we can compose this open open embedding with 
the morphism $\UU^{ns, \prime}_{1,n}\to \UU^{ns}_{1,n}$
constructed in Proposition \ref{U-prime-prop}.

Now we recall (see \cite[Lem.\ 3.1]{Smyth-I}) that every $m$-stable curve has a fundamental decomposition 
$$C=E\cup R_1\cup\ldots\cup R_k,$$
where $E$ is the minimal {\it elliptic subcurve}  i.e., the connected subcurve of arithmetic genus $1$ without disconnecting nodes, and each $R_i$ is a connected nodal curve of arithmetic genus $0$ meeting $E$ in a unique point
such that $E\cap R_i$ is a node of $C$ (and $R_i\cap R_j=\emptyset$ for $i\neq j$). 

We claim that there is at least one marked point $p_i$ on $E$. Indeed, otherwise
the $m$-stability of $C$ implies that $k>m\ge \frac{n-1}{2}$, i.e., $k\ge \frac{n+1}{2}$. 
But each $R_i$ contains at least two marked points
(due to the last condition in the definition of $m$-stability),
so the total number of marked points is $\ge 2k>n$, which is a contradiction.

Let $p_i\in E$. Then one has
$$H^1(C,\OO_C(p_i))=H^1(E,\OO_E(p_i))=0.$$
Indeed, the vanishing of $H^1(E,\OO(p_i))$ can be deduced from the classification of the possible minimal elliptic subcurves
(see \cite[Lem.\ 3.3]{Smyth-I}): $E$ is either a smooth elliptic curve, or an irreducible nodal curve, or a wheel of projective
lines, or an elliptic $l$-fold curve (which includes the cuspidal curve for $l=1$). 
Hence, we have $h^1(C,\OO(p_1+\ldots+p_n))=0$, as required.
\ed

\begin{rem}\label{Smyth-rem}
Recall that for each $i=1,\ldots,n$ we have an open subset 
$\wt{\UU}_{1,n}(i)=\pi^{-1}(U_i)\sub\wt{\UU}^{ns}_{1,n}$ consisting of
$(C,p_\bullet,v_\bullet)$ such that $h^1(p_i)=0$. The intersection 
$\cap_{i=1}^n\wt{\UU}_{1,n}(i)$ corresponds to the curves $(C,p_\bullet,v_\bullet)$ such that
$h^1(p_i)=0$ for each $i$. As shown in \cite[Prop.\ 3.3.1]{P-krich}, there is a natural projection
$$\cap_{i=1}^n\wt{\UU}_{1,n}(i)\to \wt{\UU}^{sns}_{1,n},$$ 
which is a $\G_m^{n-1}$-torsor, where the space $\wt{\UU}^{sns}_{1,n}$ classifies $(C,p_\bullet,\om)$,
where $C$ is of arithmetic genus $1$, $h^1(p_i)=0$ for each $i$ and $\om$ is a nonzero global section of
the dualizing sheaf on $C$ (and $\OO(p_1+\ldots+p_n)$ is ample). The space $\wt{\UU}^{sns}_{1,n}$
was studied in \cite{LP}, where we showed in particular that considering the GIT-quotients of 
$\wt{\UU}^{sns}_{1,n}$ by the $\G_m$-action rescaling $\om$ we recover Smyth's moduli space 
of $(n-1)$-stable curves $\ov{\MM}_{1,n}(n-1)$.
Proposition \ref{Smyth-moduli-prop} suggests that more generally, for 
$m\ge \frac{n-1}{2}$, there should exist natural morphisms from
the Smyth's moduli spaces $\ov{\MM}_{1,n}(m)$ to some GIT quotients of $\wt{\UU}_{1,n}^{ns}$ by $\G_m^n$.
We will explore this elsewhere.
\end{rem}

\section{$A_\infty$-moduli}\label{ainf-sec}

\subsection{Nice quotients}\label{nice-quotients-sec}

Below we work with schemes and group schemes over a fixed base scheme $S$.

\begin{defi}\label{nice-quot-def}
Let $G$ be a group scheme, $X$ be a $G$-scheme. We say that a $G$-invariant morphism $\pi:X\to Q$ is a {\it nice quotient}
for the $G$-action on $X$ if locally over $S$ (in Zariski topology) there exists a section $\si:Q\to X$ of $\pi$ and a morphism
$\rho:X\to G$, such that 
\begin{equation}\label{nice-quot-eq}
x=\rho(x)\si(\pi(x)) \ \ \text{ and } \ \ \rho\circ\si=1.
\end{equation}
We say that $\pi$ is a {\it strict nice quotient} if $\rho$ and $\si$ can be defined globally over $S$.
\end{defi}

In the case when $S$ is a point we obtain precisely the situation of \cite[Def.\ 4.2.2]{P-ainf}, where we called $\si(Q)$ a
{\it nice section} for the action of $G$ on $X$.

Note that a nice quotient is automatically a categorical quotient (in the category of $S$-schemes).
Indeed, let $f:X\to Z$ be a $G$-invariant morphism. Then $f(x)=f(\si(\pi(x)))$, so $f$ is a composition of 
$f\circ\si:Q\to Z$ with $\pi$.
This implies that the existence of a nice quotient is a local quesion in $S$. Namely, if $X_i\to Q_i$ are nice quotients for
$X_i=p^{-1}(U_i)$, where $(U_i)$ is an open covering of $S$, $p:X\to S$ is a projection, then we can glue them into
a global morphism $\pi:X\to Q$. 

\begin{rem} If $\pi:X\to Q$ is a nice quotient for the $G$-action on $X$ then 
$\pi$ is a universal geometric quotient (see \cite{MF}). Indeed, any base change of $\pi$ is still a nice quotient.
The following properties are clear: $\pi$ is surjective, $U\sub Q$ is open if and only if $\pi^{-1}(U)$ is open,
geometric fibers are precisely the orbits of geometric points.
Finally, we claim that $\OO_Q$ coincides with $G$-invariants in $\pi_*\OO_X$.
Indeed, given a $G$-invariant function $f$ on $\pi^{-1}(U)$ then $f(x)=f(\si(\pi(x))$, so it descends to the function
$f\circ\si$ on $U$.
\end{rem}


Let us consider the topology on the category $\Sch_S$ of $S$-schemes, such that open coverings
of $p:T\to S$ are pull-backs under $p$ of Zariski open coverings of $S$. We call this {\it $S$-Zariski topology}.
Let us consider the following presheaf of sets on $\Sch_S$:
$$T\mapsto X(T)/G(T).$$

\begin{lem} Let $\pi:X\to Q$ is a nice quotient for the $G$-action then 
the sheafification of the above presheaf with respect to the $S$-Zariski topology 
is naturally isomorphic to the functor represented by $Q$. 
Thus, a $T$-point of $Q$ can be represented by a collection of $V_i$-points of $X$, where $V_i=f^{-1}(U_i)$
for some open covering $(U_i)$ of $S$, such that for any $i$, $j$, the corresponding $V_{ij}$-points of $X$,
where $V_{ij}=f^{-1}(V_i\cap V_j)$, differ by $G(V_{ij})$-action. 
\end{lem}

\Pf . We have a natural morphism from $X(T)/G(T)$ to the sheaf represented by $Q$, which becomes an
isomorphism over an open affine covering of $S$ (due to the existence of a decomposition \eqref{nice-quot-eq}).
This immediately implies the assertion.
\ed

We have the following relative analog of \cite[Lem.\ 4.2.3]{P-ainf}. 

\begin{lem}\label{nice-quotients-main-lem}
Let $G$ be a group scheme over $S$ acting on a scheme $X$ over $S$. Assume that
$G$ fits into an exact sequence of group schemes
$$1\to H\to G\to G'\to 1$$
and that the projection $G\to G'$ admits a section $s:G'\to G$ which is a morphism of schemes
(not necessarily compatible with the group structures).
Suppose we have a scheme $Y$ with an action of $G'$ and a morphism $f:X\to X'$
compatible with the $G$-action via the homomorphism $G\to G'$.
Assume that there exists a nice quotient $\pi_H:X\to Q_H$ for the $H$-action on $X$ and 
a nice quotient $\pi':X'\to Q'$ 
for the $G'$-action on $X'$. Finally, assume that the following condition holds: for
any $S$-scheme $T$ and any points $x\in X(T)$, $g\in G(T)$ such that $f(gx)=f(x)$ there exists
an open covering $T=\cup T_i$ and a point $h\in H(T_i)$ for each $i$, such that $gx=hx$.
Then there exists a nice quotient for the $G$-action on $X$. The same assertion holds for strict nice quotients.
\end{lem}

\Pf . It is enough to prove the assertion for strict nice quotients.
Without loss of generality we can assume that the section $s:G'\to G$
satisfies $s(1)=1$. 
By assumption, we have 
sections $\si_H:Q_H\to X$ and $\si':Q'\to X'$ and the corresponding maps
$\rho_H:X\to H$ and $\rho':X'\to G'$ satisfying \eqref{nice-quot-eq}. 
Let us define morphisms $\rho_f:X\to G$ and $\pi_f:X\to X$ by
$$\rho_f=s\circ\rho'\circ f, \ \ \pi_f(x)=\rho_f(x)^{-1}x.$$
One immediately checks that 
$$f\circ \pi_f=\si'\circ\pi'\circ f.$$
In particular, $\pi_f(x)\in f^{-1}(\si'(Q'))$.
Let us set $\wt{Q}=f^{-1}(\si'(Q'))\sub X$. 
Note that for $x\in \wt{Q}$ we have
$$\rho_f(x)=s(\rho'(f(x))=s(1)=1,$$
since $\rho'|_{\si'(Q')}=1$. 
Hence, for $x\in \wt{Q}$ we have
$\pi_f(x)=x$. 
Now we set 
$$Q=\si_H^{-1}(\wt{Q})\sub Q_H,$$ 
and define the maps $\pi:X\to Q$ and $\rho:X\to G$ required for the definition of a nice quotient by
$$\pi=\pi_H\circ \pi_f,$$
$$\rho(x)=\rho_f(x)\rho_H(\pi_f(x)).$$
Note that $\pi$ is well-defined. Indeed, we need to show that
$(\si_H\pi_H\pi_f)(x)\in \wt{Q}$. But $\pi_f(x)\in\wt{Q}$, so this follows from the identity
$$(\si_H\pi_H\pi_f)(x)=\rho_H(\pi_f(x))^{-1}\pi_f(x)$$
and the fact that $\wt{Q}$ is preserved by the action of $H$.
As in \cite[Lem.\ 4.2.3]{P-ainf}, we check that our data defines a strict nice quotient for the $G$-action on $X$
(where the section of $\pi$ is provided by $\si_H|_Q$).
\ed

\subsection{General $A_\infty$-moduli}\label{gen-ainf-moduli-sec}

For a graded sheaf $\FF$ of locally free $\OO$-modules over a scheme $S$ we denote by $CH^{s+t}(\FF/S)_t$ the sheaf
of homomorphisms of $\OO$-modules $\FF^{\ot s}\to \FF$ of degree $t$. We have a natural
notion of an $A_n$-structure (resp., $A_\infty$-structure) on $\FF$, given by a collection of global sections
$$m=(m_1,\ldots,m_n)\in H^0(S,CH^2(\FF/S)_1\times\ldots\times CH^2(\FF/S)_{2-n})$$ 
(resp., $m=(m_1,m_2,\ldots)$ with $m_n\in CH^2(\FF/S)_{2-n}$),
satisfying the standard $A_\infty$-identities involving only $m_1,\ldots,m_n$ (resp., all $A_\infty$-identities).
Note that in the case when $2$ is invertible on $S$
the identities for $A_n$-structures can be written as $\sum_{i=1}^r [m_i,m_{r+1-i}]=0,$
for $r=1,\ldots,n$, where $[\cdot,\cdot]$ is the Gerstenhaber bracket.

We are interested in {\it minimal} $A_n$-structures (resp., $A_\infty$-structure), i.e., those with $m_1=0$.
The action of the group of gauge transformations on the set of $A_n$-structures also immediately generalizes
to the relative context: we have a sheaf of groups $\fG$ over $S$, where
an element of $\fG(U)$ is a collection of sections
$$f=(f_1=\id,f_2,\ldots)\in H^0(U, CH^1(\FF/S)_{-1}\times CH^1(\FF/S)_{-2}\times\ldots),$$
with the product rule obtained by interpreting $f$ as a coalgebra automorphism of the bar-coalgebra of $\FF$
(see \cite[Def.\ 4.1.3]{P-ainf}). 
We use the notation $\fG[2,n-1]:=\fG/\fG_{\ge n}$, 
introduced in \cite[Sec.\ 4.2]{P-ainf}, for the quotient of $\fG$
acting on the set of minimal $A_n$-structures on $\FF$. We denote the projection $\fG\to \fG[2,n-1]$ by
$u\mapsto u_{\le n-1}$.

\begin{rem}
The above definition of an $A_n$-algebra over a scheme is a bit naive. 
A more flexible notion should involve defining $m_i$'s only over an open
covering $U_i$ of $S$, and the gluing should be given by a collection of higher homotopies defined on intersections 
$U_{i_1}\cap\ldots\cap U_{i_r}$. We do not need the most general definition since we only aim at constructing the usual
space as a moduli space of $A_\infty$-structures (in good situations), not an $\infty$-stack. Even at this
level we will need a certain gluing procedure, but a much simpler one.
\end{rem}

Now let us fix a scheme $S$ and a sheaf $\EE$ of associative $\OO_S$-algebras over $S$. We assume also that $\EE$ is locally free of finite rank over $\OO_S$. Roughly speaking, our goal is to classify gauge equivalence classes of
minimal $A_\infty$-algebras that extend the associative algebras $\EE_s:=\EE|_s$, where $s$ is a point of $S$.

To begin with, for every $n\ge 2$ we define the functor $\AA_n=\AA_{n,\EE}$ (resp., $\AA_\infty=\AA_{\infty,\EE}$)
on the site of $S$-schemes, which associates with $f:T\to S$ the set of minimal $A_n$-structures (resp., $A_\infty$-structures) extending the sheaf of associative $\OO_T$-algebras $f^*\EE$.
This functor is clearly represented by an affine scheme $\AA_n(\EE)$ over $S$. Namely,
$A_n(\EE)$ is the closed subscheme in the total space of the vector bundle
$CH^2(\EE/S)_{-1}\oplus\ldots\oplus CH^2(\EE/S)_{2-n}$ given by the $A_\infty$-equations.
We have a natural projection
\begin{equation}\label{An-projection-eq}
\pi_n:\AA_n\to \AA_{n-1}: m\mapsto m_{\le n-1}.
\end{equation}

Next, we have the sheaf of groups $\fG$ of gauge transformations acting on each functor $\AA_n$ through the quotient
$\fG[2,n-1]$, and we make the following
defintion.

\begin{defi}\label{ainf-moduli-pre-def}
Let $\wt{\MM}_n$ denote the  quotient-functor associating to an $S$-scheme $f:T\to S$ the set $\AA_n(T)/\fG[2,n-1](T)$
of gauge equivalence classes of minimal $A_n$-structures on $f^*\EE$. 
We denote by $\wt{\MM}_\infty$ the similar quotient-functor for gauge equivalence classes of minimal
$A_\infty$-structures. 
\end{defi}

Note that the sheaf of groups $\fG[2,n-1]$ is representable by a unipotent affine group scheme over $S$ 
which we still denote as $\fG[2,n-1]$. Note also that the projection
$\fG\to\fG[2,n-1]$ admits a section (not compatible with the group structures) and so is universally surjective.
However, the quotient-functor $\wt{\MM}_n$ is not necessarily representable. 

\begin{lem}\label{subbun-lem} 
(i) Assume that $S$ is reduced and Noetherian, and let $(V^\bullet,d)$ be a bounded below
complex of vector bundles over $S$ such that
$H^i(V^\bullet_s)=0$ for $i<p$.
Then for $i<p$, the image $\im(d^i)$ of the differential $d^i:V^i\to V^{i+1}$
is a subbundle of $V^{i+1}$. 

\noindent
(ii) If in addition $S$ is affine then there exist decompositions
$V^i=B^i\oplus C^i$, for $i\le p$, such that for $i<p$ one has $d^i(C^i)=B^{i+1}$ and the map
$$d^i|_{C^i}: C^i\to B^{i+1}$$
is an isomorphism. In this situation for any $f:T\to S$ the complex $H^0(T,f^*V^\bullet)$ is exact
in degrees $<p$.
\end{lem}

\Pf . (i) It is enough to prove that $\im(d^{p-1})$ is a subbundle in $V^p$.
Without loss of generality we can assume that $V^i=0$ for $i<0$ and $p>0$. Then the map
$V^0_x\to V^1_x$ is injective for every $x\in S$, so $d:V^0\to V^1$ is the embedding of a subbundle.
Hence, $V^1/d(V^0)$ is a vector bundle, and we can replace our complex with 
$$0\to V^1/d(V^0)\to V^2\to\ldots $$
and iterate the same argument.

\noindent
(ii) The first assertion follows from part (i): we set $B^i:=\im(d^{i-1})$ and let $C^i$ be the image of any splitting
of the projection $V^i\to V^i/B^i$ (which exists since $S$ is affine). These decompositions carry over to the complex
$H^0(T,f^*V^\bullet)$, which implies its exactness in degrees $<p$.
\ed

We denote the Hochschild differential $[m_2,?]$ on $CH^*(\EE/S)$ by $\de$ and its graded components as
$$\de^i_t:CH^i(\EE/S)_t\to CH^{i+1}(\EE/S)_t.$$ 

\begin{lem}\label{gauge-eq-lem} 
Let $S$ be a reduced Noetherian affine scheme. 

\noindent
(i) Assume that $HH^i(\EE_s)_{-j}=0$ for $i\le 1$ and $j=1,\ldots,d-2$, for a fixed $d\ge 2$.
The for $f:T\to S$, let $m$ and $m'$ be a pair of minimal $A_n$-structures on $f^*\EE$ for some $n\ge d$, such that 
$m$ is gauge equivalent to $m'$ (over $T$) and
$m_{\le d}=m'_{\le d}$. Then there exists a gauge equivalence $u$ over $T$, such that
$u_{\le d-1}=\id$ and $m'=u\cdot m$.

\noindent
(ii) Assume that $HH^i(\EE_s)_{<0}=0$ for $i\le 1$. Then 
the natural map
$$\wt{\MM}_\infty(T)\to \liminv_n \wt{\MM}_n(T)$$
is an isomorphism for every $S$-scheme $T$. 
\end{lem}

\Pf . (i) The proof is similar to that of \cite[Lem.\ 4.1.6]{P-ainf}. By assumption, we have a gauge equivalence 
$\wt{u}\in\fG(T)$ such that $\wt{u}m=m'$. Then $\wt{u}_{\le d}$ sends $m_{\le d}$ to itself.
Now Lemma \ref{subbun-lem}(ii), applied to the Hochschild complexes $(CH^*(\EE)_{-j},\de)$, implies that
$HH^1(H^0(T,f^*\EE)/\OO(T))_{-j}=0$ for $j=1,\ldots,d-2$. Hence, arguing as in \cite[Lem.\ 4.1.6]{P-ainf},
we can correct $\wt{u}$ to a gauge equivalence $u\in \fG(T)$ such that 
$u\cdot m=m'$ and $u_{\le d-1}=\id$.

\noindent
(ii) The proof is identical to the argument in the proof of Corollary \cite[Cor.\ 4.2.5]{P-ainf}, with part (i) replacing the reference
to \cite[Lem.\ 4.1.6]{P-ainf}.
\ed

\begin{defi}\label{ainf-moduli-def}
Let us denote by $\MM_n$ (resp., $\MM_\infty$) 
the sheafification of the functor $\wt{\MM}_n$ (resp., $\wt{\MM}_\infty$) with respect to the $S$-Zariski topology.
\end{defi}

The following theorem, which is a relative version of \cite[Thm.\ 4.2.4]{P-ainf}, shows that under certain
vanishing assumptions on the Hochschild cohomology, the functor $\MM_n$ is representable by an affine $S$-scheme.
Note that these assumptions are slightly stronger than in \cite[Thm.\ 4.2.4]{P-ainf}.

\begin{thm}\label{gen-ainf-thm} 
Assume that $S$ is reduced and Noetherian, and that $HH^i(\EE_s)_{-j}=0$ for $i\le 1$ and $1\le j\le n-3$.
for every point $s\in S$. 

\noindent
(i) There exists a nice quotient $\AA_n(\EE)/\fG[2,n-1]$ for the action of $\fG[2,n-1]$ on $\AA_n(\EE)$.
This quotient $\AA_n(\EE)/\fG[2,n-1]$, which is affine of finite type over $S$, represents the functor $\MM_n$.
If in addition $S$ is affine then there exists a strict nice quotient $\AA_n(\EE)/\fG[2,n-1]$, and the natural map
of functors $\wt{\MM}_n\to \MM_n$ is an isomorphism.

\noindent
(ii) Assume $HH^i(\EE_s)_{-j}=0$ for $i\le 1$ and $j\ge 1$. Then
the scheme $\liminv_n \MM_n$, affine over $S$, represents the functor $\MM_\infty$. 
In the case when $S$ is affine,
the natural map $\wt{\MM}_\infty\to\MM_\infty$ is an isomorphism.

\noindent
(iii) Assume $HH^i(\EE_s)_{-j}=0$ for $i\le 1$, $j\ge 1$, and in addition $HH^2(\EE_s)_{-j}=0$ for $j>n-2$. Then the morphism
$\MM_\infty\to\MM_n$ is a closed embedding. If in addition $HH^3(\EE_s)_{-j}=0$ for $j>n-2$ then
$\MM_\infty\to \MM_n$ is an isomorphism.
\end{thm}

\Pf . (i) It is enough to prove the existence of a strict nice
quotient for the $\fG[2,n-1]$-action on $\AA_n=\AA_n(\EE)$ in the case when $S$
is affine. Indeed, then it would follow that $\wt{\MM}_n$ is represented by this quotient, and hence, the map
$\wt{\MM}_n\to\MM_n$ is an isomorphism.

Similarly to the proof of \cite[Thm.\ 4.2.4]{P-ainf} the existence of a strict nice quotient is proved by the induction on $n$,
using Lemma \ref{nice-quotients-main-lem}. 
Assume that $n>2$ and we already have a section
$S_{n-1}$ for the $\fG[2,n-2]$-action on $\AA_{n-1}$. We have an exact sequence of sheaves of groups over $S$,
$$0\to CH^1(\EE)_{2-n}\to \fG[2,n-1]\to \fG[2,n-2]\to 0.$$
We want to find a section for the $CH^1(\EE)_{2-n}$-action on $\AA_n$. By
Lemma \ref{subbun-lem}(ii), there exists 
a complement $\KK_{2-n}\sub CH^2_{2-n}$
to the subbundle $\im\de^1_{2-n}$.
Let $\AA'_n$ denote the closed subset of $\AA_n$ given by the condition $m_n\in \KK_{2-n}$.
Since the action of $x\in CH^1(\EE)_{2-n}$ on $(m_2,\ldots,m_n)\in \AA_n$ changes $m_n$ to $m_n+\de^1(x)$
and does not change $(m_2,\ldots,m_{n-1})$, we see that $\AA'_n$ is a section for the $CH^1(\EE)_{2-n}$-action
on $\AA_n$. Now we can apply Lemma \ref{nice-quotients-main-lem} to the projection \eqref{An-projection-eq}
and the compatible actions of $\fG[2,n-1]\to \fG[2,n-2]$.
Note that to apply this Lemma we need to check
that the intersection of an $\fG[2,n-1]$-orbit with a fiber of $\pi_n$ is a $CH^1(\EE)_{2-n}$-orbit.
But this follows from Lemma \ref{gauge-eq-lem}(i).
Thus, we deduce that $\AA'_n\cap \pi_n^{-1}(S_{n-1})$ is a section for the $\fG[2,n-1]$-action on $\AA_n$.

\noindent
(ii)  First, assume that $S$ is affine. Then, combining part (i) with Lemma \ref{gauge-eq-lem}(ii),
we derive that the functor $\wt{\MM}_\infty$ is represented by the scheme $\liminv_n \MM_n$, affine over $S$.
Hence, in this case the map $\wt{\MM}_\infty\to\MM_\infty$ is an isomorphism.
Thus, in the case of general $S$ the map of sheaves $\MM_\infty\to \liminv_n \MM_n$ becomes an isomorphism over
an affine open covering of $S$, hence, it is an isomorphism.

\noindent
(iii) We can assume $S$ to be affine. Then this is proved similarly to \cite[Cor.\ 4.2.6]{P-ainf}.
\ed


\subsection{$A_\infty$-structures associated to curves}\label{ainf-curve-sec}

We want to study the relative moduli space of $A_\infty$-structures on 
the family of associative algebras $(E_W)$ over the Grassmannian $G(n-g,n)$
(see \eqref{E-W-eq}).

First, let us describe more precisely the corresponding sheaf of $\OO$-algebras $\EE_{g,n}$ over $G(n-g,n)$.

Let $R$ be a commutative ring, and 
let $J_0$ be the two-sided ideal in the path algebra $R[Q_n]$ of $Q_n$ generated by the elements
$$A_iB_iA_i, B_iA_iB_i, A_iB_j,$$
where $i\neq j$.
Given an $R$-submodule $W\sub R^n$ such that $R^n/W$ is locally free of rank $r$,
we define $J_W\sub R[Q_n]$ as  the ideal generated by $J_0$ together with the additional relations
$$\sum x_i B_iA_i=0 \ \ \text{ for every } \sum x_i e_i\in W.$$
The corresponding quotient algebra
$$E_W=R[Q_n]/J_W$$ 
is projective as an $R$-module. In the case when $R^n/W$ is a free $R$-module, $E_W$ is also free over $R$, 
with the basis given by the elements $(A_i), (B_i), (A_iB_i)$ and some $g$ linear combinations 
$\sum x_i B_iA_i$ such that $\sum x_ie_i$ project to a basis of $R^n/W$. 

For an $n$-tuple of invertible elements $\la=(\la_1,\ldots,\la_n)=(R^*)^n$ we have a natural isomorphism
$$E_W\to E_{\la\cdot W}: A_i\mapsto A_i, B_i\mapsto \la_iB_i,$$
where the transformation $W\mapsto \la\cdot W$ is induced by the rescaling of the basis $e_1,\ldots,e_n$ of $R^n$.

Applying the above construction to the tautological subbundle over an affine covering of the Grassmannian
$G(n-g,n)$ and gluing, we obtain a $\G_m^n$-equivariant sheaf $\EE_{g,n}$ of $\OO$-algebras over $G(n-g,n)$,
where $\G_m^n$ acts naturally on $G(n-g,n)$ (as diagonal matrices).

\begin{defi}\label{M-infty-defi} 
The moduli functor $\MM_\infty$ on $G(n-g,n)$-schemes associates with $f:S\to G(n-g,n)$, 
a minimal $A_\infty$-structure on $f^*\EE_{g,n}$, given locally over an open covering of $S$ (with compatibility up
to a gauge equivalence on intersections), viewed up to a gauge equivalence.
\end{defi}


Similarly to \cite[Sec.\ 3]{P-ainf} we have a natural morphism of functors
\begin{equation}\label{curves-to-ainf-map}
\wt{\UU}^{ns}_{g,n}\to \MM_\infty.
\end{equation}
Namely, to a family of curves in $\wt{\UU}^{ns}_{g,n}$ over an affine base we associate the corresponding
minimal $A_\infty$-structure, defined using some relative formal parameters along the marked points and
the dg-model and homotopies described in \cite[Sec.\ 3]{P-ainf}. 
Note that if we choose different formal parameters with the same underlying tangent vectors,
then we will get the same dg-algebra but with different homotopies, hence the corresponding $A_\infty$-structures will be gauge 
equivalent.
Thus, the map \eqref{curves-to-ainf-map} is well defined. 
Futhermore, it is compatible with the $\G_m^n$-actions, where the $\G_m^n$-action on $\MM_\infty$ is induced by
the rescalings
$$A_i\mapsto A_i, \ B_i\mapsto \la_i B_i,$$
for $(\la_1,\ldots,\la_n)\in\G_m^n$.

Next, we want to prove that $\MM_\infty$ is represented by an affine scheme of finite type over $G(n-g,n)$.
For this we want to apply the criterion of Theorem \ref{gen-ainf-thm}, which requires some information about
the Hochschild cohomology of the algebras $E_W$. As in \cite{P-ainf}, we will get this information geometrically
by identifying $HH^*(E_W)$ with the Hochschild cohomology of the corresponding special curve. 

\begin{lem}\label{formal-lem} Let $C_W=(C(\ov{h}),p_\bullet,v_\bullet)\in\wt{\UU}^{ns}_{g,n}$
be the special curve corresponding to $W\in G(n-g,n)$ (see Proposition \ref{special-curves-prop}). Then
the natural $A_\infty$-structure on $\Ext(G,G)$ for $G$ given by \eqref{G-generator-eq} is trivial (up to a gauge equivalence),
 and
hence, we have an isomorphism
$$HH^*(C_W)\simeq HH^*(E_W),$$
where $W=\pi(C,p_1,\ldots,p_n)$. The second grading on $HH^*(E_W)$ corresponds to the weights of the $\G_m$-action,
coming from the natural $\G_m$-action on $C_W$.
\end{lem}

\Pf . This is similar to \cite[Prop.\ 4.4.1]{P-ainf}.
\ed

\begin{lem}\label{tangent-weights-lem}
Let $C$ be a reduced projective curve over a field $k$ with a $\G_m$-action, which is
the union of irreducible components $C_i$, $i=1,\ldots,n$, joined in a single point $q$. Assume
that $C\setminus\{q\}$ is smooth and that each normalization map $\wt{C}_i\to C_i$ is a bijection,
with $\wt{C}_i\simeq\P^1$.
Assume also that the action of $\G_m$ on the Zariski tangent space at $q$ has negative weights.
Then 

\noindent
(i) the action of $\G_m$ on $H^1(C,\OO_C)$ has positive weights.

\noindent
(ii) Assume in addition that $C=C_W$ for some subspace $W\sub k^n$, where $W=0$ if $n=1$.
Assume also that $\cha(k)\neq 2$, and if $n=1$ then also $\cha(k)\neq 3$.
Then $H^1(C,\TT)=0$, and the action of $\G_m$ on $H^0(C,\TT)$ has weights $0$ and $1$.

\noindent
(iii) Keep the assumptions of (ii). Let $p_i\in C_i\setminus \{q\}$ be the unique $\G_m$-invariant point, and let
$D=\sum_i p_i$. Then one has 
$$H^0(C,\TT(-D))=H^0(C,\TT)^{\G_m}, \ \ H^0(C,\TT(-2D))=0.$$
Also the natural map $H^0(C,\TT(nD))\to H^0(C,\TT(nD)|_D)$ is surjective for $n\ge 0$.

\noindent
(iv) For $W=k^n$, with $n\ge 2$, the assertions of (ii) and (iii) hold without any restrictions on the characteristic of $k$.
\end{lem}

\Pf . (i) Let $V=C\setminus\{q\}$. We can choose a coordinate $x_i$ on an affine part of $\wt{C}_i\simeq\P^1$ containing $q$
such that $x_i(q)=0$ and $x_i$ has positive weight $w_i$ with respect to the $\G_m$-action. 
Let $U$ be an affine neighborhood of $q$ obtained by deleting on each $C_i$ the point where $x_i$ has a pole.
We can calculate $H^1(C,\OO_C)$ as the quotient of $\OO(U\setminus\{q\})$
by $\OO(V)+\OO(U)$. Since every $x_i^n$ with $n\le 0$ extends to a regular function on $V$, we see that
$H^1(C,\OO_C)$ is spanned by positive powers of $x_i$'s, so $\G_m$ has only positive weights on it.

\noindent
(ii) The case $n=1$, $W=0$ corresponds to the cuspidal curve, for which the assertions of (ii) and (iii) are known
(see e.g., \cite[Lem.\ 4.4.2]{P-ainf}). So we assume $n\ge 2$.
We use the coordinates $x_i$ on affine parts of the normalizations $\wt{C}_i$ from Definition \ref{sp-curve-def}.
The space $H^0(C,\TT_C)$ embeds into the space of vector fields on $V\simeq \sqcup_{i=1}^n (C_i\setminus \{q\})$, which are spanned by $x_i^m\del_{x_i}$ with $m\le 2$. 
Exactly as in the proof of \cite[Lem.\ 4.4.2(ii)]{P-ainf} we check that if a vector field $v=(P_i(x_i,x_i^{-1})\del_{x_i})$
on $U\setminus \{q\}$ extends to a derivation of $\OO(U)$ then $P_i\in x_ik[x_i]$ for every $i$
(this uses the assumption $\cha(k)\neq 2$). Thus, if $v$ extends to a global section of $\TT_C$ then
each $P_i$ is a linear combination of $x_i$ and $x_i^2$, which implies that the weights of $\G_m$ on $H^0(C,\TT_C)$
are $0$ and $1$.
Similarly, we see that if each $P_i\in x_i^2k[x_i]$ then $v$ extends to a derivation of $\OO(U)$.
Thus, $H^0(U,\TT)$ and $H^0(V,\TT)$ span $H^0(U\setminus\{q\},\TT)$, which gives the vanishing of $H^1(C,\TT_C)$.

\noindent
(iii) A vector field on $U\setminus \{q\}$ has zero (resp., double zero) along $D$ iff each $P_i\in x_ik[x_i^{-1}]$
(resp., $P_i\in k[x_i^{-1}]$). Together with calculations of (ii) this immediately implies our assertions about
$H^0(C,\TT(-D))$ and $H^0(C,\TT(-2D))$. Next,
similarly to (ii) we can represent sections of $H^0(C,\TT_C(nD))$ as vector fields
$v=(P_i(x_i)\del_{x_i})$ on $U\setminus\{q\}$ with $\deg(P_i)\le n+2$, and the last assertion follows from the fact that
$v$ extends to a regular derivation of $\OO(U)$ whenever $P_i\in x_i^2k[x_i]$.

\noindent
(iv) We can argue as in the proof of \cite[Lem.\ 4.4.2(ii)]{P-ainf}, 
with $x_i^2$ replaced by $x_i$ (complemented by zeros in all other places), to show
that vector fields on $V$ that extend to $U$ are precisely $v=(P_i\del_{x_i})$ with $P_i\in x_ik[x_i]$. The rest of the
proof is the same as in (ii) and (iii).
\ed

\begin{rem}\label{special-W-rem}
One can check that the assertions of Lemma \ref{tangent-weights-lem} hold for $C=C_W$ without any restrictions
on the characteristic of $k$, provided $n\ge g+2$ and $W$ is not contained in any of the coordinate hyperplanes
$k^{n-1}\sub k^n$.
\end{rem}

\begin{cor}\label{HH1-cor} 
For any subspace $W\sub k^n$, where $W=0$ if $n=1$, and 
$k$ is a field of characteristic $\neq 2$ (resp., $\neq 2,3$ if $n=1$), one has
$$HH^0(E_W)_{<0}=HH^1(E_W)_{<0}=0.$$
The same result holds for $W=k^n$, $n\ge 2$, with no restrictions on the characteristic.
\end{cor}

\Pf . By Lemma \ref{formal-lem}, we have $HH^1(E_W)\simeq HH^1(C_W)$, where $C_W$ is the corresponding
special curve, and the second grading is induced by the $\G_m$-action on $C_W$. Now $HH^0(C_W)=H^0(C_W,\OO)$
lives in degree $0$. For $HH^1$ we use the
exact sequence
$$0\to H^1(C_W,\OO)\to HH^1(C_W)\to H^0(C_W,\TT)\to 0$$
(see \cite[Sec.\ 4.1.3]{LPer}). Now the assertion follows
from Lemma \ref{tangent-weights-lem}(i)(ii).
\ed

\begin{prop}\label{ainf-mod-prop} 
Let us work over $\Z[1/2]$ if $n\ge 2$, or over $\Z[1/6]$ if $n=1$, or over $\Z$ if $g=0$. 
Assume that either $n\ge 2$ or $g=1$.
Then the functor $\MM_\infty$ of $A_\infty$-structures (up to a gauge equivalence) on the family $(E_W)$ 
is represented by an affine scheme of finite type over $G(n-g,n)$.
\end{prop}

\Pf . Due to Corollary \ref{HH1-cor}, the criterion of Theorem \ref{gen-ainf-thm}(ii) implies that $\MM_\infty$ (resp., $\MM_n$)
is represented by an affine scheme (resp., of finite type) over $G(n-g,n)$.
Next, we note that by Lemma \ref{formal-lem}, $HH^i(E_W)$ is 
finite-dimensional for every $i$. Hence, by Theorem \ref{gen-ainf-thm}(iii), we derive that $\MM_\infty\simeq\MM_n$ for
sufficiently large $n$, so it is of finite type over $G(n-g,n)$.
\ed

For a scheme $S$ over a field $k$ we denote by $L_S$ the cotangent complex of $S$ over $k$.

\begin{lem}\label{cot-com-lem} 
Assume that either $W=k^n$ and $n\ge 2$, or $k$ has characteristic $\neq 2$ (resp., $\neq 2,3$ if $n=1$).
Let $C=C_W$ be a special curve over $k$, where $W=0$ if $n=1$, and let $D=p_1+\ldots+p_n$, $U=C\setminus D$.
Then the natural morphism 
\begin{equation}\label{Ext1-L-C-map}
\Ext^1(L_C,\OO_C(-2D))\to \Ext^1(L_C,\OO_C(-D))
\end{equation}
is surjective, while the natural morphism
$$\Ext^1(L_C,\OO_C(-D))\to \Ext^1(L_U,\OO_U)$$
is an isomorphisms.
The natural morphism
$$\Ext^2(L_C,\OO_C(-2D))\to \Ext^2(L_U,\OO_U)$$
is an isomorphism.
\end{lem}

\Pf . The proof is almost the same as that of \cite[Lem.\ 4.5.4]{P-ainf}, using Lemma \ref{tangent-weights-lem}. 
The difference is that in our case the map
$$H^0(C,\TT(-D))\to H^0(C,\TT(-D)|_D)$$
is not necessarily surjective, so we cannot assert that the map \eqref{Ext1-L-C-map} is an isomorphism,
only that it is surjective.
\ed

Now we can use the results of \cite{P-ainf} to compare the deformation theory of a special curve $C_W$
with that of $E_W$, viewed as an $A_\infty$-algebra.

We refer to \cite{Man} for the basic deformation theory and for some terminology used below.
Let us fix a field $k$ and consider the category $\Art(k)$ of local Artinian algebras with the residue field $k$. 
Given a curve $(C,p_1,\ldots,p_n,v_1,\ldots,v_n)$ with smooth distinct marked points and the nonzero
tangent vectors at them, we have the corresponding deformation functor
$$\Def(C,p_\bullet,v_\bullet):\Art(k)\to \Sets$$
associating with $R$ the set of isomorphism classes of flat proper families of curves $\pi_R:C_R\to\Spec(R)$ with sections
$p^R_1,\ldots,p^R_n$, and trivializations of the relative tangent bundle along them, such that the induced
data over $\Spec(k)\sub\Spec(R)$ is $(C,p_\bullet,v_\bullet)$.

On the other hand, for any finite-dimensional minimal $A_\infty$-algebra $E$ we have the deformation functor 
$$\Def(E):\Art(k)\to \Sets$$
of extended\footnote{In an extended gauge transformation $(f_1,f_2,\ldots)$ the map $f_1$ is only required
to be invertible, not necessarily equal to $\id$, see \cite[Def.\ 4.1.3]{P-ainf}.}
gauge equivalence classes of minimal $A_\infty$-algebras $E_R$ over $R$, reducing to $E$ over $k$.
Let also, for a fixed $n-g$-dimensional subspace $W\sub k^n$, 
$$\wt{\Def}(E_W):\Art(k)\to \Sets$$
be the functor associating with $R$ the set of pairs $(W_R,m_\bullet)$, where $W_R$ is an $R$-point of $G(n-g,n)$, 
reducing to $W$ over $k$, and $m_\bullet$ is a minimal $A_\infty$-structure on $E_{W_R}$, reducing to the trivial
$A_\infty$-structure on $E_W$, viewed up to a gauge equivalence.
Note that the functor $\wt{\Def}(E_W)$ is prorepresented by the formal completion of the scheme
$\MM_\infty$ at the point corresponding to the trivial $A_\infty$-structure on $E_W$.
We have a natural forgetful morphism 
$$\wt{\Def}(E_W)\to \Def(E_W).$$ 

\begin{lem}
The tangent space to the functor $\wt{\Def}(E_W)$ can be identified with 
$$HH^2(E_W)_{<0}\oplus T_W G(n-g,n).$$
There is a complete obstruction theory for this functor with values in $HH^3(E_W)_{<0}$.
\end{lem}

\Pf . The tangent space classifies pairs $(f,m_\bullet)$, where $f:\Spec(k[t]/(t^2))\to G(n-g,n)$ is a morphism 
sending the closed point to $W$, and $m_\bullet$ is a minimal $A_\infty$-structure on $f^*\EE$, extending the given $m_2$,
reducing to the trivial one modulo $(t)$,
up to a gauge equivalence. Then $f$ corresponds to a tangent vector in $T_W G(n-g,n)$, while the class of
$(m_3,m_4,\ldots)$ is an element in $HH^3(E_W)_{<0}$ (see e.g., \cite[Lem.\ 4.5.2]{P-ainf}).
The obstruction theory is obtained from the usual obstruction theory for $A_\infty$-structures 
(see e.g., \cite[Lem.\ 4.5.2]{P-ainf}) using the fact that $G(n-g,n)$ is smooth.
\ed

For each special curve $(C_W,p_\bullet,v_\bullet)\in\wt{\UU}^{ns}_{g,n}$ 
corresponding to a subspace $W\in G(n-g,n)(k)$, where $k$ is a field,
the morphism \eqref{curves-to-ainf-map} induces a morphism of deformation functors
\begin{equation}\label{def-functors-map}
\Def(C_W)\to \wt{\Def}(E_W).
\end{equation}

\begin{prop}\label{deform-prop} Assume that either $n\ge 2$ and $g=0$, or
$n\ge 2$ and the characteristic of $k$ is $\neq 2$, or $n=g=1$ 
and the characteristic of $k$ is $\neq 2, 3$. Then
the morphism \eqref{def-functors-map} is an isomorphism.
\end{prop}

\Pf . Let $U_W$ be the affine curve $C_W\setminus D$, where $D=p_1+\ldots+p_n$.
Let also $\CC$ be  the (non-full) subcategory in the $A_\infty$-enhancement of the derived category
of $\Qcoh(C_W)$
with the objects $(\OO_C,\OO_{p_1},\ldots,\OO_{p_n},\OO_{U_W})$, all morphisms (including $\Ext^*$) between $\OO_C,(\OO_{p_i})$, and all morphisms from $\OO_C$ to $\OO_{U_W}$ and from $\OO_{U_W}$ to $\OO_{U_W}$ (but we don't consider any other morphisms between these objects).
As in the proof of \cite[Prop.\ 4.5.4]{P-ainf},
we consider additional functors $\Def(U_W)$, $\Def_{nc}(U_W)$ and $\Def(\CC)$ 
(where they are denoted $F_{U_W}$, $F_{U_W,nc}$ and $F_{\CC}$), corresponding to deformations
of $U_W$, of $\OO(U_W)$ as an associative algebra, and of $\CC$ as an $A_\infty$-category.
These functors fit into a commutative diagram
\begin{equation}\label{def-fun-diagram}
\begin{diagram}
\Def(U_W)&\lTo{}& \Def(C_W) &\rTo{} & \wt{\Def}(E_W)\\
\dTo{}&&\dTo{}&&\dTo{}\\
\Def_{nc}(U_W)&\lTo{}& \Def(\CC) &\rTo{}& \Def(E_W)
\end{diagram}
\end{equation}

\noindent
{\bf Step 1}. The map $\Def(\CC)\to \Def(E_W)$ is \'etale. This is proved in the same way as in Step 1 of
\cite[Prop.\ 4.5.4]{P-ainf}.

\noindent
{\bf Step 2}. The map $\Def(C_W)\to \Def(U_W)$ is smooth, while the map $\Def(U_W)\to \Def_{nc}(U_W)$ is \'etale.

First, we observe that the maps on tangent spaces induced by these maps are
$$\Ext^1(L_{C_W},\OO(-2D))\to \Ext^1(L_{U_W},\OO_{U_W})\to HH^2(U_W),$$
the first of which is surjective by Lemma \ref{cot-com-lem}, while the second is an isomorphism by 
\cite[Lem.\ 4.4.6]{P-ainf}.
Similarly the maps of obstruction spaces are
$$\Ext^2(L_{C_W},\OO(-2D))\to \Ext^2(L_{U_W},\OO_{U_W})\to HH^3(U_W),$$
of which the first is an isomorphism by Lemma \ref{cot-com-lem}, while the second is injective by \cite[Lem.\ 4.4.6]{P-ainf}.
Hence, the maps $\Def(C_W)\to \Def(U_W)$ and $\Def(U_W)\to \Def_{nc}(U_W)$ are smooth 
and the second is \'etale (see \cite[Prop.\ 2.17]{Man}).

\noindent
{\bf Step 3}. The map $\Def(C_W)\to \Def(\CC)$ (resp., $\Def(\CC)\to \Def_{nc}(U_W)$) induces a surjection 
(resp., isomorphism) on tangent spaces.

Indeed, Step 2, together with the commutativity of diagram \eqref{def-fun-diagram}, 
implies that $\Def(\CC)\to \Def_{nc}(U_W)$ induces a surjection on tangent spaces. But 
$$HH^2(U_W)\simeq HH^2(C_W)\simeq HH^2(E_W),$$ 
so the dimensions of tangent spaces are the same. Hence, $\Def(\CC)\to \Def_{nc}(U_W)$ induces an isomorphism
on tangent spaces.

It follows that the maps induced on tangent spaces by $\Def(C_W)\to \Def(\CC)$ and by $\Def(C_W)\to \Def_{nc}(U_W)$
are isomorphic, so the required surjectivity follows from Step 2.

\noindent
{\bf Step 4}. The map $\Def(C_W)\to \wt{\Def}(E_W)$ (resp., $\wt{\Def}(E_W)\to \Def(E_W)$) induces an isomorphism
(resp., surjection) on tangent spaces.

Note that by Steps 1 and 3, we know that the map $\Def(C_W)\to \Def(E_W)$ induces a surjection on tangent spaces.
Hence, the same is true for $\wt{\Def}(E_W)\to \Def(E_W)$.
We claim that there is a commutative diagram with exact rows
\begin{equation}\label{tangent-maps-diagram}
\begin{diagram}
k^n&\rTo{\a}& \Ext^1(L_{C_W},\OO(-2D)) &\rTo{\b}& HH^2(\CC)  &\rTo{} 0\\
\dTo{=}&&\dTo{\ga}&&\dTo{\ga'}\\
k^n&\rTo{\a'}& HH^2(E_W)_{<0}\oplus T_W G(n-g,n)&\rTo{\b'}& HH^2(E_W)_{\le 0} &\rTo{} 0
\end{diagram}
\end{equation}
where the arrow $\a$ (resp., $\a'$) is induced by the $\G_m^n$-action on the functor $\Def(C_W)$ (resp.,
$\wt{\Def}(E_W)$), while the right commutative square is induced by the right commutative square in \eqref{def-fun-diagram}
(flipped about the diagonal).
Note that we already know that $\ga'$ is an isomorphism and $\b'$ is surjective. 
To see the exactness of the top row we observe that by Steps 1 and 2, the map $\b$ can be identified with the morphism
$$\Ext^1(L_{C_W},\OO(-2D))\to \Ext^1(L_{C_W},\OO(-D))\rTo{\sim} \Ext^1(L_{U_W},\OO_{U_W}),$$
where the second arrow is an isomorphism by Lemma \ref{cot-com-lem}. Hence, its kernel is the image of the coboundary
map $H^0(C_W,\TT(-D)|_D)\to\Ext^1(L_{C_W},\OO(-2D))$, which can be identified with $\a$.
The exactness of the bottom row in \eqref{tangent-maps-diagram}
would follow from the exactness in the middle of the sequence
$$k^n\to T_W G(n-g,n)\to HH^2(E_W)_0\to 0,$$
where the second arrow is the tangent map to the map $W\to E_W$, and the first arrow corresponds to the $\G_m^n$-action
on $G(n-g,n)$. But this follows from the observation 
that a $k[t]/(t^2)$-point of $G(n-g,n)$, $\WW$, can be recovered from the isomorphism class of the corresponding
algebra $E_\WW$ up to a $\G_m^n$-action.

Note that diagram \eqref{tangent-maps-diagram}, together with the fact that $\ga'$ is an isomorphism, immediately implies
that $\ga$ is surjective. It remains to prove that the restriction of $\ga$ to $\im(\a)$ is injective. To this end we use the 
fact that each point $C_W\in \wt{\UU}^{ns}_{g,n}$ lies in the section $\si(G(n-g,n))$ of the projection to $G(n-g,n)$,
and that the $\G_m^n$-orbit of $C_W$ still lies in $\si(G(n-g,n))$. Hence, the tangent space to this orbit maps injectively
to $T_W G(n-g,n)$, which implies our assertion.

\noindent
{\bf Step 5}. The morphisms $\Def(C_W)\to \Def_{nc}(U_W)$, $\Def(C_W)\to \Def(\CC)$ and $\Def(C_W)\to \Def(E_W)$ are smooth, and the morphism $\Def(C_W)\to \wt{\Def}(E_W)$ is an isomorphism.

The first morphism is equal to the composition
$$\Def(C_W)\to \Def(U_W)\to \Def_{nc}(U_W),$$
where both arrows are smooth by Step 2. But it is also equal to the composition
$$\Def(C_W)\to \Def(\CC)\to \Def_{nc}(U_W)$$
By Step 3, the first arrow induces a surjection on tangent spaces. Hence, by \cite[Lem.\ 4.5.3]{P-ainf}, the morphism
$\Def(C_W)\to \Def(\CC)$ is smooth. Using Step 1 again we deduce that $\Def(C_W)\to \Def(E_W)$ is smooth.
Thus, the composition
$$\Def(C_W)\to \wt{\Def}(E_W)\to \Def(E_W)$$
is smooth, and the first arrow induces an isomorphism of tangent spaces. Hence, by \cite[Lem.\ 4.5.3]{P-ainf},
the morphism $\Def(C_W)\to \wt{\Def}(E_W)$ is smooth, and hence, \'etale.
But the functor $\wt{\Def}(E_W)$ is homogeneous, since it is prorepresented (by the formal completion of the scheme
$\MM_\infty$), so it is an isomorphism by \cite[Cor.\ 2.11]{Man}.
\ed

\begin{thm}\label{ainf-thm} 
Assume that either $n\ge g\ge 1$, $n\ge 2$ and the base is $\Spec(\Z[1/2])$, or
$n=g=1$ and the base is $\Spec(\Z[1/6])$, or $g=0$, $n\ge 2$ and the base is $\Spec(\Z)$.
Then the morphism \eqref{curves-to-ainf-map} is an isomorphism.
\end{thm}

\Pf . We know that both schemes are affine of finite type over $G(n-g,n)$
(by Theorem A and Proposition \ref{ainf-mod-prop}), and that the morphism
\eqref{curves-to-ainf-map} is compatible with $\G_m$-action. 
Furthermore, the $\G_m$-invariant loci of each scheme provide a section of the projection to $G(n-g,n)$.
Thus, locally over over $G(n-g,n)$ our morphism corresponds to a homomorphism $f:A\to B$
of non-negatively graded algebras such that $f_0:A_0\to B_0$ is an isomorphism. Furthermore, by
Proposition \ref{deform-prop}, for every 
point of $\Spec(A_0)\simeq\Spec(B_0)$, the map $f$ induces an isomorphism of deformation functors.
Hence, applying Lemma \ref{graded-lem} below we deduce that $f$ is an isomorphism.
\ed

\begin{lem}\label{graded-lem} 
Let $f:A\to B$ be a morphism of degree zero of non-negatively graded algebras such that the induced map
$A_0\to B_0$ is an isomorphism. Assume that $A_0$ is Noetherian, $A$ and $B$ are finitely generated as
algebras over $A_0\simeq B_0$, and for every maximal ideal $\mg\sub A_0$ the map $f$ induces
an isomorphism $\hat{A}\to \hat{B}$ of the completions with respect to the maximal ideals $\mg+A_{>0}$ and
$\mg+B_{>0}$, respectively. Then $f$ is an isomorphism.
\end{lem}

\Pf . It is enough to prove
that $f$ induces an isomorphism $A/A_{>0}^N\to B/B_{>0}^N$ for each $N>0$. 
Note that $A/A_{>0}^N$ (resp., $B/B_{>0}^N$) is a finitely generated module over $A_0$ (resp., $B_0$).
Note that for any maximal ideal $\mg\sub A_0\simeq B_0$,
the $(\mg+A_{>0})$-adic topology on $A/A_{>0}^N$ is equivalent to the $\mg$-adic topology, and similarly on
$B/B_{>0}^N$.
Thus, we have a morphism 
$$A/A_{>0}^N\to B/B_{>0}^N$$
of finitely generated $A_0$-modules, inducing an isomorphism of $\mg$-adic
completions of localizations at every maximal ideal $\mg\sub A_0$.
Since $A_0$ is Noetherian, such a morphism is an isomorphism.
\ed

\begin{rem} 
For $g\ge 1$ let us define the open subset $U\sub G(n-g,n)$ to be the complement
to the union of the images of $n$ embeddings $G(n-g,n-1)\sub G(n-g,n)$ associated with the coordinate 
hyperplanes $k^{n-1}\hra k^n$. It is easy to see that the preimage $\pi^{-1}(U)\sub \wt{\UU}^{ns}_{g,n}$ parametrizes
$(C,p_\bullet,v_\bullet)$ such that $H^1(C,\OO(D-p_i))=0$ for every $i$ (where $D=p_1+\ldots+p_n$).
Using Remark \ref{special-W-rem} one can see that the analog of Proposition \ref{ainf-mod-prop} gives a
relative moduli of $A_\infty$-structures on the family $(E_W)$ over $U$, when working over $\Spec(\Z)$,
provided $n\ge g+2$. Similarly, Proposition \ref{deform-prop} holds without any restrictions on the characteristic, for
$W\in U$ and $n\ge g+2$. This suggests that for $n\ge g+2$, working over $\Spec(\Z)$, 
one could still show that the morphism $\pi^{-1}(U)\to U$ is affine of finite type. Then the analog of Theorem \ref{ainf-thm}
would give an isomorphism of $\pi^{-1}(U)$ with the corresponding relative moduli of $A_\infty$-structures over $U$.
\end{rem}

\end{document}